\documentclass[preprint, 12pt]{elsarticle}
\usepackage{amssymb}
\usepackage{amsmath} 
\usepackage{bm}
\usepackage{mathrsfs}
\usepackage{soul}
\usepackage{array}
\usepackage{longtable}
\usepackage{caption}
\usepackage{subcaption}
\usepackage{float}
\usepackage{tikz}
\usepackage{parskip}
\usepackage[a4paper, left=2.5cm, right=2.5cm, top=3cm, bottom=3cm]{geometry}
\usepackage[hidelinks]{hyperref} 

\usepackage{setspace}
\usepackage{cleveref}
\usepackage{tcolorbox}
\newtcolorbox{mybox}[1]{fonttitle=\bfseries,title=#1}

\usepackage{xcolor}
\usepackage{soul}

\graphicspath{{./figures/}}

\newcommand{\leff}{$l_{\mathrm{eff}}$}
\newcommand{\heff}{$h_{\mathrm{eff}}$}

\newcommand{\fmax}{$F_{\mathrm{max}}$}

\newcommand{\crossA}{\mathcal{A}}

\renewcommand{\hl}{}

\begin{document}

\begin{frontmatter}

\title{On the failure of beam-like topologically interlocked structures}
\author[inst1]{Ioannis Koureas}
\author[inst1]{Mohit Pundir}
\author[inst1]{Shai Feldfogel}
\author[inst1]{David S. Kammer\corref{cor1}}
\cortext[cor1]{Corresponding author}
\ead{dkammer@ethz.ch}
\affiliation[inst1]{organization={Institute for Building Materials}, 
            addressline={ETH Zurich},
            city={Zurich},
            postcode={8093},
            country={Switzerland}}

\begin{abstract}
Topologically interlocked structures are architectured by fitting together blocks that are constrained geometrically through contact and friction by their neighboring blocks. As long as the frictional strength is nowhere exceeded, the blocks stick against each other, allowing for large rotations. Once the interfacial stresses exceed the frictional strength, relative sliding between the blocks alters the structure's mechanical response. Improving the structural performance, precisely the strength and the toughness, has been one of the main focal points in the literature. However, many fundamental questions regarding the role and effect of the interface mechanisms (stick and slip) and rotation of the blocks have not been addressed yet. Here, we carry out a parametric analysis to understand the effect of Young's modulus, friction coefficient and geometry of the blocks on the dominance of the stick or slip governed mechanism. We combine analytical and computational tools to analyze the failure mechanisms and the response capacities of beam-like topologically interlocked structures. This is achieved using a finite element method coupled with a penalty-based approach for enforcing contact constraints along interfaces. We show that the combination of the structure's height and the friction coefficient controls whether the failure mechanism is slip-governed or stick-governed. Furthermore, we demonstrate that the sticking mechanism across all interfaces along with the rotation of the blocks dictates a saturation level to the mechanical performance of a given structure irrespective of geometric and material properties. This provides a theoretical upper bound for the structural response of topologically interlocked structures and establishes a theoretical benchmark of achievable performance.
\end{abstract}



\end{frontmatter}

\section{Introduction}
\label{sec:introduction}

\begin{figure}[H]
    \begin{center}
    \begin{tikzpicture}
    \node[inner sep = 0] at (0, 0) {\includegraphics[width=0.9\linewidth]{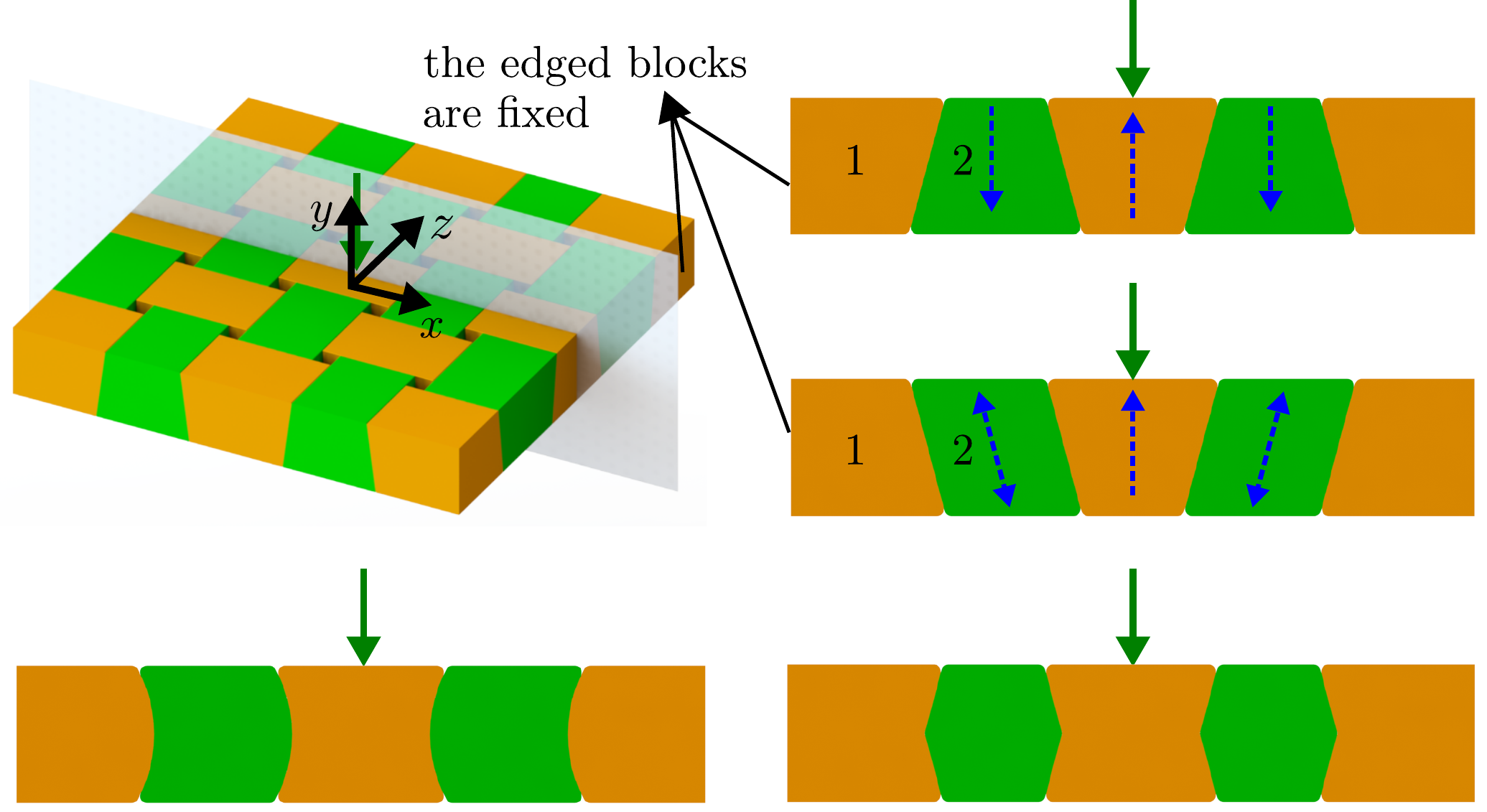}};
    \node[inner sep = 0] at (-6.8, 3.5) {(a)};
    \node[inner sep = 0] at (0.7, 3.5) {(b)};
    \node[inner sep = 0] at (0.7, 0.8) {(c)};
    \node[inner sep = 0] at (-6.8, -2.1) {(d)};
    \node[inner sep = 0] at (0.7, -2.1) {(e)};
    \end{tikzpicture}
    \caption{Rationale for examined configuration: (a) A topologically interlocked structure (TIS) assembly with planar-faced blocks~\citep{Mirkhalaf2019}. (b) Central cross-section in the x-y plane and (c) our configuration obtained from the original x-y cross-section with the inclination between blocks 1 and 2 flipped. Schematic representation of beam-like structure with (d) curved and (e) kinked surfaces.}
    \label{fig:tim-structure}
    \end{center}
\end{figure}
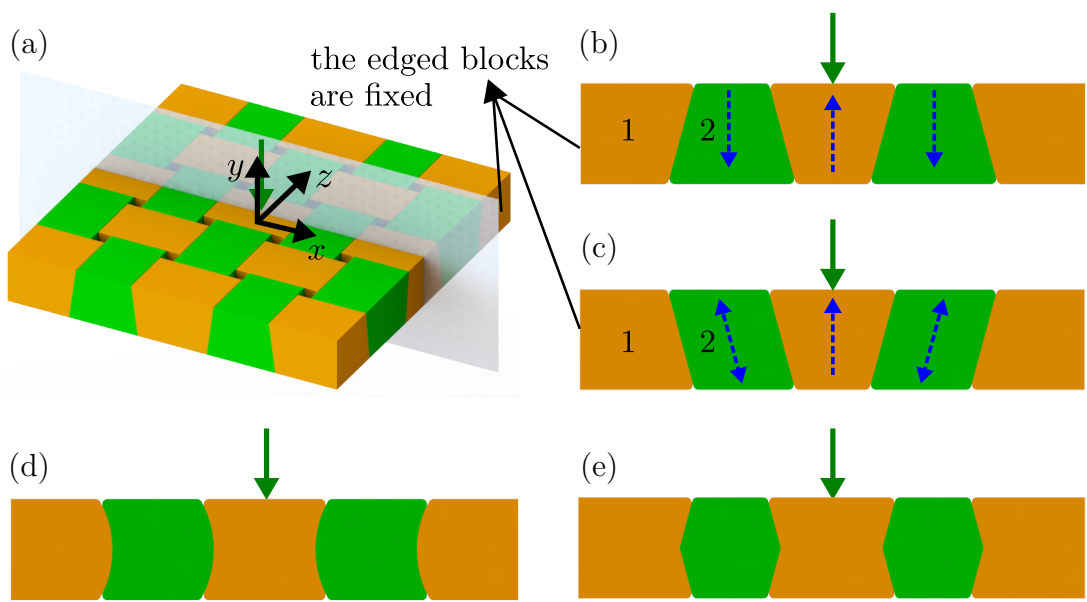

Topologically interlocked structures (TIS) are assemblies of building blocks that hold together due to the blocks' unique interlocking shapes \citep{Dyskin2001a, Dyskin2001, Chen2008, Djumas2016, Yin2019, Dyskin2003a, Schaare2008}. The unbonded nature of the blocks means that TIS rely on contact and frictional interactions between blocks for structural integrity~\mbox{(\autoref{fig:tim-structure}a)}. TIS enjoy unique structural properties, including high toughness against failure for structures made from brittle material~\citep{Djumas2016, Mirkhalaf2014, Mirkhalaf2018} and structural integrity despite partial failure (\textit{e.g.}, missing blocks in plate-like TIS)~\citep{Dyskin2003}. \hl{However, TIS have not yet found widespread application in engineering because their highly complex, non-linear behavior and failure are not yet fully understood, hindering the ability to design them safely. Better understanding the behavior and failure of TIS is therefore required and it is at the focus of this study.}

The main mechanisms governing the behavior and failure of TIS are local slip and stick combined with large block rotations~\citep{Mirkhalaf2018, Khandelwal2015, Mirkhalaf2019}. These mechanisms take more or less prominent roles in the structural response depending on the TIS design and material properties. \hl{However, it remains unclear how these properties determine which mechanism will likely be the more dominant one in a given TIS configuration and how they affect the global structural response in terms of load-carrying capacity, loading energy, global stiffness and ultimate deflection.}

\begin{figure}[H]
    \begin{center}
    \begin{tikzpicture}
    \node[inner sep = 0] at (0, 0) {\includegraphics[width=0.5\linewidth]{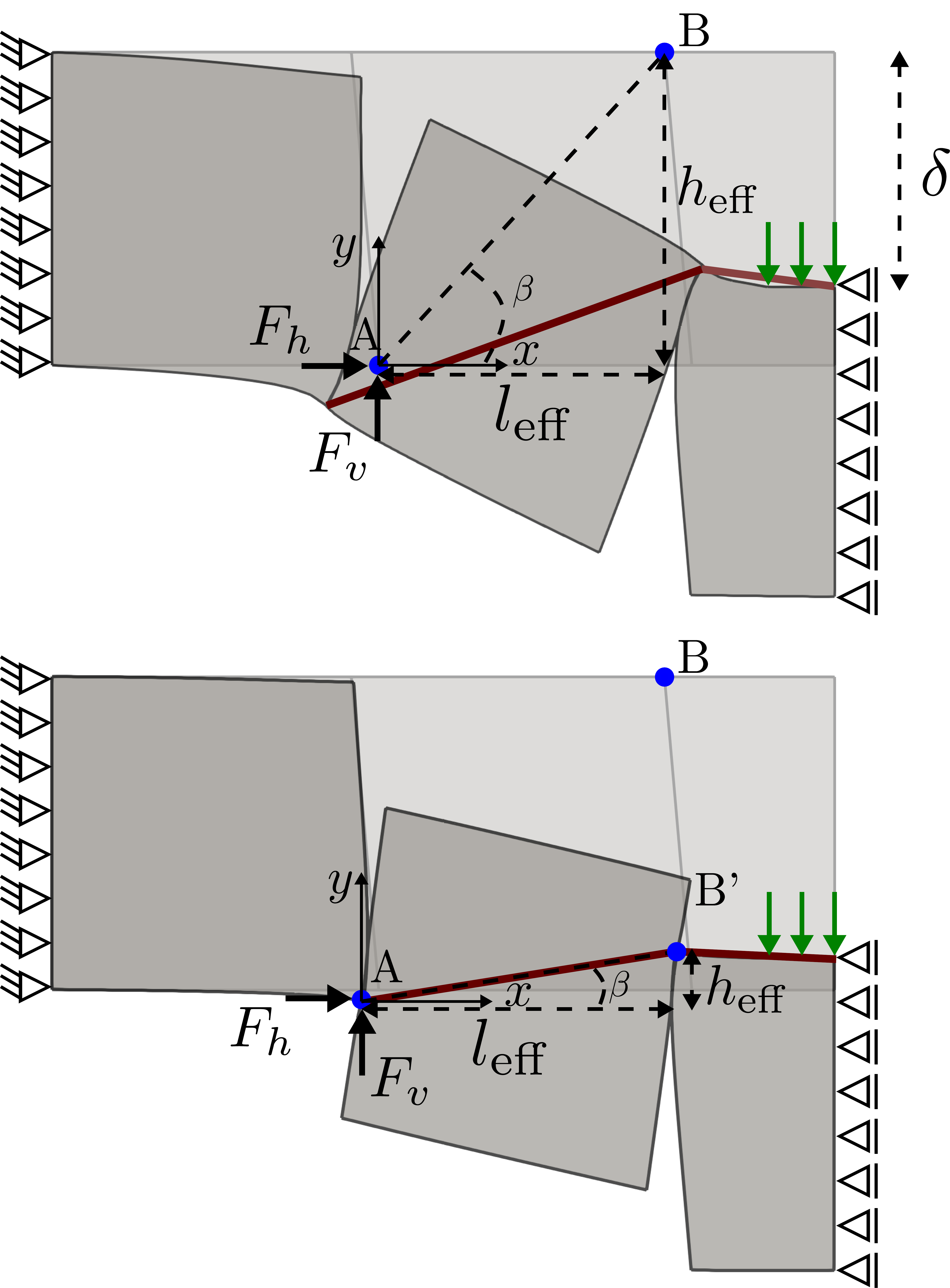}};
    \node[inner sep = 0] at (-4.3, 5.5) {(a)};
    \node[inner sep = 0] at (-4.3, 0.2) {(b)};
    \end{tikzpicture}
    \caption{Equivalent truss model (dark red color) in a topologically interlocked structure (TIS) for (a) a stick and (b) a slip scenario. The effective height $h_{\mathrm{eff}}$ and length $l_{\mathrm{eff}}$, as well as the displacement $\delta$ are shown for the deformed structures.}
    \label{fig:truss-model}
    \end{center}
\end{figure}

In the relatively simple case where the response is entirely stick-governed, the behavior can be described by analytical models as in~\citep{Khandelwal2015, Khandelwal2014, Khandelwal2012}. Khandelwal \textit{et~al.} modeled TIS as a truss and derived analytical expressions for the horizontal and vertical reaction forces $F_h$ and $F_v$ at the local pivoting point A~(\mbox{\autoref{fig:truss-model}a}), which goes as follows:
\begin{equation} \label{eq:truss model}
    F_v = F_h \left(\frac{h_{\mathrm{eff}} - \delta}{l_{\mathrm{eff}}} \right) = E \crossA \frac{l_{\mathrm{eff}}}{\sqrt{l_{\mathrm{eff}}^2 + h_{\mathrm{eff}}^2}} \left( \frac{\sqrt{l_{\mathrm{eff}}^2 + h_{\mathrm{eff}}^2}}{\sqrt{l_{\mathrm{eff}}^2 + (h_{\mathrm{eff}} - \delta)^2}} - 1\right) \left( \frac{h_{\mathrm{eff}} - \delta}{l_{\mathrm{eff}}} \right)
\end{equation}
where $E$ is the Young's modulus of the material, \heff{} and \leff{} are the effective height and length respectively~\mbox{(\autoref{fig:truss-model})}, $\crossA$ is the cross-sectional area of the truss model and $\delta$ the applied displacement if no slip occurred so far. \autoref{eq:truss model} clearly demonstrates that in a sticking situation the load-carrying capacity (\textit{i.e.}, $2F_v$) scales linearly with $E$ and is strongly dependent on \heff{} and \leff{}.  Although~\autoref{eq:truss model} does not apply when slip occurs, we can still qualitatively explain how local slip reduces global stiffness by examining a specific slipped configuration. Specifically, when slip occurs, the initial pivoting point B moves to B'~\mbox{(\autoref{fig:truss-model}b)}, reducing \heff{}, $F_h$ and hence the load-carrying capacity $2F_v$. \hl{However, these observations remain qualitative, particularly regarding the possible co-existence of stick and slip. Specifically, they do not account for the occurrence of local slips and the consequent evolution of slip-governed failure mechanism. To capture a slip-governed response, which is the most common one in practice and the more challenging one to model, computational methods are required.}

The finite element method (FEM) has been shown to capture and quantify the experimentally observed failure and the load-displacement curve in beam-like TIS~\citep{Dalaq2019, Mahoney2022}. Dalaq \textit{et~al.}~\citep{Dalaq2019} showed that the failure mode depends on the number of blocks, friction coefficient and the shape of the interfaces, specifically that curved interfaces (a similar configuration is shown in Figure~\mbox{\ref{fig:tim-structure}d}) promote sliding of the blocks and delay hinging (i.e., stick and rotation)~\citep{Dalaq2020}. \hl{However, Dalaq \textit{et~al.}~\mbox{\citep{Dalaq2020}} focused on a specific material and did not investigate the effects of Young's modulus $E$ on the failure mechanisms and the structural capacity. In addition, the effect of the assembly's height $h$, which governs the bending stiffness of TIS assemblies, see Zakeri \textit{et~al.}~\mbox{\citep{Zakeri2021}}, was not considered.}

\hl{In summary, the effects of $h$, $E$ and $\mu$ on the type of failure, \textit{i.e.}, a slip-governed or stick-governed one and on the associated structural capacities of TIS have hitherto not been addressed in the literature and remain only partially understood. Here, we aim to clarify and quantify these effects. Specifically, we aim to better understand how different combinations of $h$, $E$ and $\mu$ tend to make the response more slip-, or stick-governed and how they affect the structural capacity in terms of maximal load, loading energy, global stiffness and ultimate deflection. Towards that aim, we perform a three-way $E-\mu-h$ parametric study using FEM, based on the latter's ability to capture and quantify experimentally observed failure mechanisms, in particular the slip-governed one, in beam-like TIS \mbox{\citep{Dalaq2019, Dalaq2020}}. In contrast with previous parametric studies that were limited to the stick-regime, the main strength of the present study is that it treats both the stick and the slip regimes within a unified FEM-based framework. This allows us not only to better understand the previously unaddressed effects of $E$, $\mu$ and $h$ on the slip-governed failure, but also to better understand the conditions that control the two mechanisms.}

In the following, we discuss the choice of examined configuration and the FEM formulation (section~\ref{sec:numerical model}). In section~\ref{sec:results and discussion} we present and discuss the results of the parametric study and give an outlook on how the gained knowledge can aid the design of TIS.

\section{Numerical Model}\label{sec:numerical model}

\subsection{Examined configuration}\label{sec:examined configuration}

\hl{To address the effects of $E$, $\mu$ and $h$ on the failure and response capacity of TIS, we choose a beam-like configuration inspired by the centrally loaded plate-like TIS studied experimentally in~\mbox{\citep{Mirkhalaf2019}}. The original configuration from~\mbox{\citep{Mirkhalaf2019}} is presented in Figure~\mbox{\ref{fig:tim-structure}a}. The central cross-section of the structure in the $x-y$ plane is depicted in Figure~\mbox{\ref{fig:tim-structure}b}. The configuration we use for this study is modified compared to the cross-section of the actual 3D TIS such that the angle of inclination of the interface between blocks $1$ and $2$ is reversed, as shown in Figure~\mbox{\ref{fig:tim-structure}c}}.

\hl{The concept of simplifying 3D TIS to 2D beam-like structures with a representative cross-section of the 3D equivalent goes back to~\mbox{\citep{Khandelwal2015, Khandelwal2014, Khandelwal2012, Dalaq2019, Dalaq2020}}. This approach, well-established in structural analysis of monolithic structures, is motivated in the present TIS context by computational and methodological considerations. The computational cost of modeling slip-governed failure in TIS with FEM is always very high  and in TIS with more than a few blocks, may be prohibitive~\mbox{\citep{Dalaq2020}}. Methodologically, we assume that the effects of material properties ($E$ and $\mu$) and the structural height ($h$) on the structural response are qualitatively similar in 3D and 2D TIS, as they are in monolithic counterparts, therefore, we consider the latter. The structural action of beams is much simpler and it allows to study these effects in pure form in accordance with the objectives of the present study.}

\hl{We note that our configuration is not fully interlocked because there are directions (indicated by blue arrows in Figures~\mbox{\ref{fig:tim-structure}b~and~c}) in which the blocks are not kinematically constrained by the neighboring blocks. This scenario is due to the 3D nature of topological interlocking and it is very typical of 2D cross-section representation, like the actual cross-section in Figure~\mbox{\ref{fig:tim-structure}b} and other beam-like TIS configurations studied in the literature~\mbox{\citep{Khandelwal2015, Khandelwal2014, Khandelwal2012, Dalaq2020}}. To effectively constrain our configuration, we required that (a) the structure has structural integrity under its self-weight and (b) that it is interlocked under the examined loads. Requirement (a) motivated reversing the angle of inclination between blocks 1 and 2 compared to the actual cross-section of the TIS. Requirement (b) is met in our configuration since we only examine a downward load on the central block. By avoiding loads in the few specific directions that are not kinematically constrained (the degrees of freedom in which our configuration deviates from a strictly defined TIS), our configuration may be considered to be effectively TIS under the examined loads. In order to examine the generality of the observation, we consider two additional cases with non-planar interfaces, a curved interface~\mbox{(\autoref{fig:tim-structure}d)} and a kinked interface~\mbox{(\autoref{fig:tim-structure}e)}.}


\begin{table}
\begin{center}
\caption{Symbols and notations}
\begin{tabular}{ l l l }
    Symbol & Description & Unit \\
    \hline
    $E$ & Young's modulus & $\mathrm{N/m^2}$ \\
    $F_a$ & Resultant force at a given pivoting point of the truss model & N \\
    $F_y$ & Load-carrying capacity of the structure & N \\
    $F_{\mathrm{max}}$ & Maximum load-carrying capacity of the structure & N \\
    $h$ & Height of a block ($=$ structural depth) & m \\
    \heff{} & Effective height in a TIS & m \\
    $K$ & \textit{Normalized} global stiffness of the structure & - \\
    $l$ & Length of a block & m \\
    \leff{} & Effective length in a TIS & m \\
    $N$ & Normal force along a TIS interface & N \\
    $t$ & Thickness of a TIS & m \\
    $T$ & Tangential force along a TIS interface & N \\
    $\delta$ & Prescribed displacement at the top central surface of the middle block & m \\
    $U$ & Loading energy of the structure & Nm \\
    $\beta$ & Angle controlled by \heff{} and \leff{} & $^o$ \\
    $\varepsilon_n$ & Normal penalty parameter & $\mathrm{N/m^3}$ \\
    $\varepsilon_t$ & Tangential penalty parameter & $\mathrm{N/m^3}$ \\
    $\theta$ & Inclined angle & $^o$ \\
    $\mu$ & Friction coefficient & - \\
    $\mu_{\mathrm{sat}}$ & Saturated friction coefficient & - \\
    $\nu$ & Poisson's ratio & - \\
\end{tabular}
\label{tab:symbols-notations}
\end{center}
\end{table}

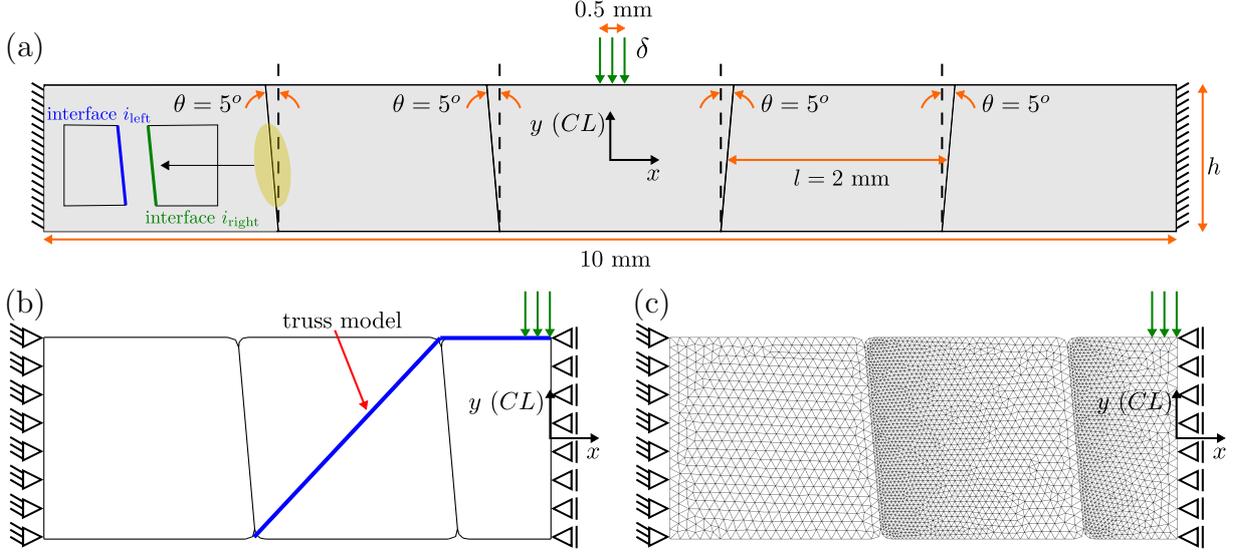
\begin{figure}[H]
    \begin{center}
    \begin{tikzpicture}
    \node[inner sep = 0] at (0, 0)
    {\includegraphics[width=1\linewidth]{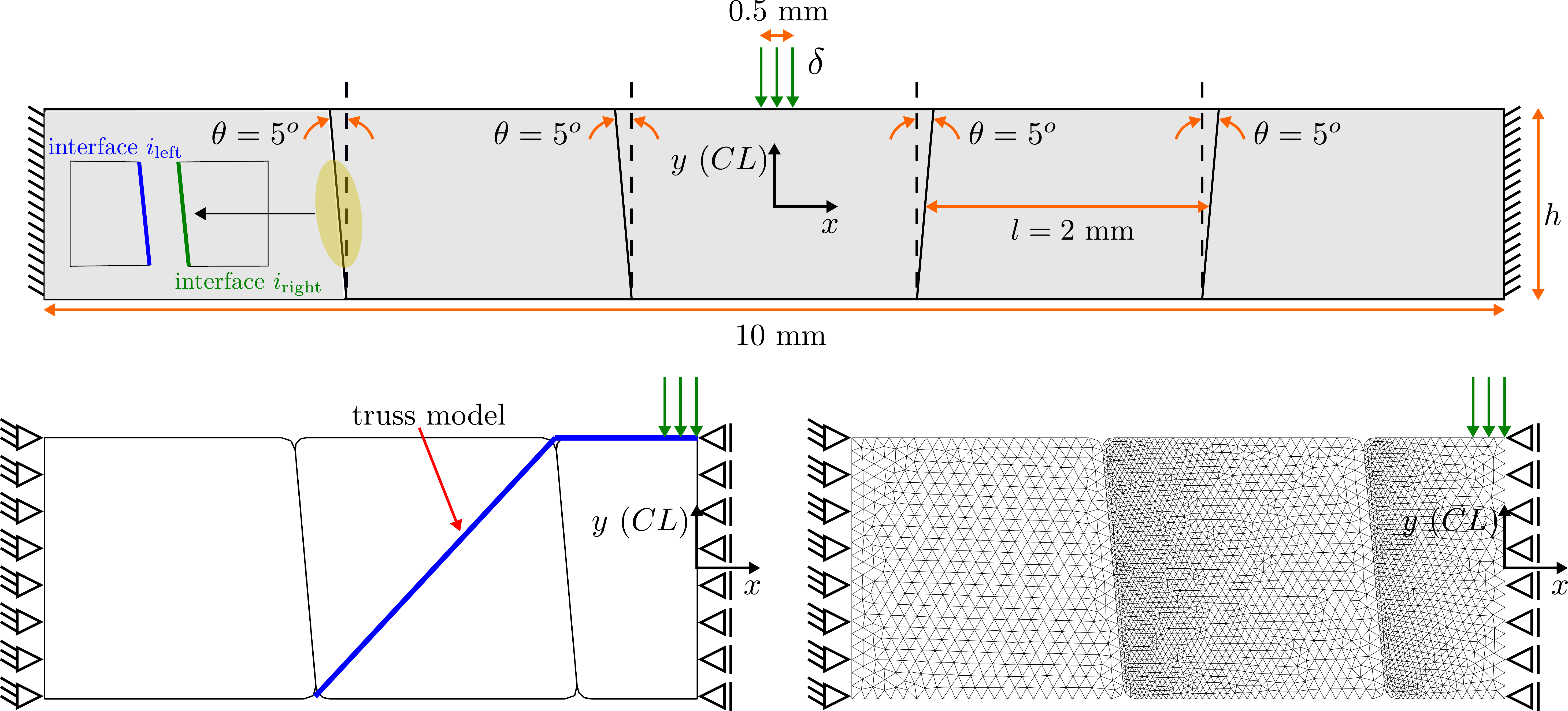}};
    \node[inner sep = 0] at (-7.8, 3) {(a)};
    \node[inner sep = 0] at (-7.8, -0.4) {(b)};
    \node[inner sep = 0] at (0.45, -0.4) {(c)};
    \end{tikzpicture}
    \caption{Schematic representation of model set-up showing (a) the geometric parameters and boundary conditions in a five-block TIS. Every structure consists of $i$ interfaces where $i = 1, 2, ..., k$ with $k$ being the total number of interfaces in a structure. Every interface consists of two sides, the left and  the right. (b) Schematic illustration of the symmetric model with respect to the center line (CL). The blue line represents the truss model. (c) Mesh and boundary conditions as used in the simulation.}
    \label{fig:tim-blocks}
    \end{center}
\end{figure}

\subsection{Numerical formulation}
\label{subsec:Finite Element formulation}

FEM is used for the numerical analyses in this study. We employ finite strain formulation to account for the large deformations and large rotations of the building blocks~\citep{Bathe1975}. Thus, considering $n$ deformable bodies $\Omega_n^{i}$, the weak formulation at load increment $i$ is described as:
\begin{equation}
\sum^{n}{\int_{\Omega^{i}_{n}}{\bm{\underline{\underline{\varepsilon}}}^i : \mathbb{C} : \delta\bm{\underline{\underline{\varepsilon}}}^i \ d\Omega^{i}_{n}}} \ + \ \sum^{n}{\int_{\Omega^{i}_{n}}{\bm{\underline{\underline{S}}}^i : \delta\bm{\underline{\underline{\eta}}} \ d\Omega^{i}_{n}}} \ = \ W_{\mathrm{ext}}^{i+1} \ - \ \sum^{n}{\int_{\Omega^{i}_{n}}{\bm{\underline{\underline{S}}}^i : \delta\bm{\underline{\underline{\varepsilon}}} \  d\Omega^{i}_{n}}}
\end{equation}
where, $\bm{\underline{\underline{S}}}$ and $\bm{\underline{\underline{\varepsilon}}}$ are the $2^{nd}$ Piola-Kirchhoff stress tensor and the linear strain tensor, respectively. $W_{ext}^{i+1}$ is the virtual work of the external forces, $\bm{\underline{\underline{\eta}}}$ represents the nonlinear incremental strain tensor and $\mathbb{C}$ the $4^{th}$ order constitutive tensor. A node-to-segment contact algorithm, with penalty-based constraints, is employed~\citep{Konyukhov2004, Konyukhov2005, Laursen2002, Wriggers2007, Yastrebov2013, Zavarise2009} to enforce contact and frictional constraints along the interfaces of $n$ deformable bodies. We use the penalty method for its computational simplicity. The virtual work $\delta W_c$ of the contact forces at the current configuration for $n$ deformable bodies that come in contact at $k$ interfaces $S_{k}$ is expressed as:
\begin{equation}\label{eq:contact}
\delta W_c = \sum^{k}{\int_{S_{k_{\mathrm{slave}}}}{(T_n \bm{n} + \bm{T}_t) \cdot (\delta \bm{u}_{k_{\mathrm{slave}}} - \delta \bm{u}_{k_{\mathrm{master}}})\ d{S_{k_{\mathrm{slave}}}}}}
\end{equation}
where $T_n$ is the traction along the normal $\bm{n}$ to the interface and $\bm{T_t}$ is the frictional traction tangential to $\bm{n}$ integrated over one of the two contact surfaces termed as slave surface. Based on the penalty approach,  $T_n=\varepsilon_n  \langle g\rangle$ is approximated as a linear function of the orthonormal gap between a slave node and the master surface. Similarly, $\bm{T}_t = \varepsilon_t (\bm{\Delta u}_t)$ is approximated as a linear function of the tangential slip distance ($\bm{\Delta u}_t$) between a slave node and the master surface, computed based on the covariant derivative approach~\citep{Konyukhov2005, Wriggers2001}. The penalty parameters ($\varepsilon_n, \varepsilon_t$) are area regularized to ensure that the computed contact forces are mesh independent~\citep{El-Abbasi2001, Zavarise2009a}.  In order to overcome the biases in choosing a slave and master surface at an interface (see~\autoref{eq:contact}), a two-pass algorithm~\citep{El-Abbasi2001, Zavarise2009a} is employed whereby at each load increment, the contact forces at a node are computed considering once a surface as a slave and then as a master. The FE code with the finite strain formulation and the node-to-segment contact algorithm, is developed as in-house code and has been validated for frictional cases (more details provided in~\mbox{\ref{app:contact-validation}}).

\hl{The examined structure is depicted in Figure~\mbox{\ref{fig:tim-blocks}a}. We consider a five-block assembly with a span length of $10$~mm and three different heights ($h = 1$~mm, $h = 1.5$~mm and $h = 2$~mm). It is fixed at its ends (i.e., $u_y(\pm 5, y) = u_x(\pm 5, y) = 0$) and it is loaded incrementally by prescribing the displacement $u_y\left(-0.25 \leq x \leq 0.25, \frac{h}{2}\right) = \delta = \frac{h}{500}$. The total force corresponding to $\delta$ is denoted by $F_y$. Using the symmetry about the $y$-axis, we model only the left half of the structure, where the symmetry boundary condition $u_x = 0$ is prescribed along $x=0$.}

Each block is characterized by its angle $\theta$, height $h$ and length $l$. The blocks are considered to be isotropic linear elastic material with Young's modulus $E$, Poisson's ratio $\nu$ and friction coefficient $\mu$. A description of the symbol notation we use is provided in~\mbox{\Cref{tab:symbols-notations}}. The material and geometrical values used for the parametric analyses are presented in~\mbox{\Cref{tab:material-parameters}}. We chose values of $E$ that cover an essential range of brittle materials and allow us to study the effect of material elasticity on the global stiffness of TIS. Furthermore, we explore a wide range of $\mu$ to understand the effect of interfacial friction of the blocks. \hl{The effects of fracture are not accounted.}

For the numerical analysis, the topologically interlocked beam~\mbox{(\autoref{fig:tim-blocks}b)} is modelled as a 2D structure under plane-strain conditions. The structure is discretized using first-order triangular elements~\mbox{(\autoref{fig:tim-blocks}c)} and the corners are rounded to avoid non-physical stress singularities. \hl{All simulations are performed under static conditions. Therefore, dynamic effects of friction and inertial effects that may be associated with the structural response are not considered.}

\begin{table}[H]
\begin{center}
\caption{Material parameters}
\renewcommand{\arraystretch}{1.1}
\begin{tabular}{ l l }
    Parameter & Value \\
    \hline
    $E~(\mathrm{GPa})$ & $1$, $2$, $3$, $10$, $20$, $30$ \\
    $h$~(mm) & 1.0, 1.5, 2.0 \\
    $l$~(mm) & 2.0 \\
    $t$~(mm) & 1.0\\
    $\theta~(^o)$ & 5 \\
    $\mu$ & 0.2, 0.4, 0.6, 0.8, 1.0, 1.2 \\
    $\nu$ & 0.2 \\
\end{tabular}
\label{tab:material-parameters}
\end{center}
\end{table}

\begin{figure}[H]
    \begin{center}
    \includegraphics[width=0.4\linewidth]{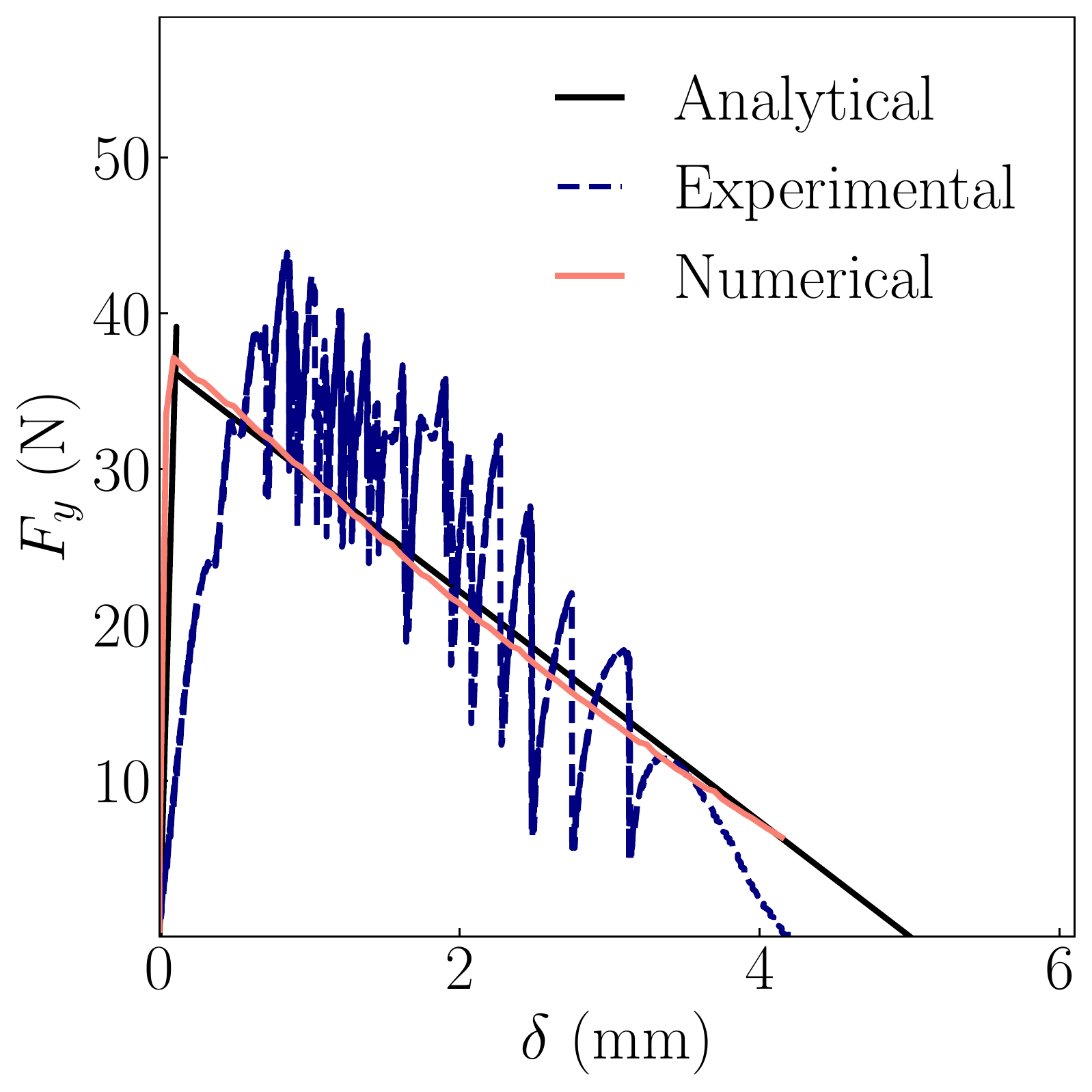}
    \caption{Load-carrying capacity $F_y$ against the prescribed displacement $\delta$ for a five-block pre-compressed structure with angle $\theta = 0^o$. The experimental and analytical results have been taken from~\citep{Dalaq2019}.}
    \label{fig:experimental-validation}
    \end{center}
\end{figure}

To validate our FE formulation in the context of TIS, we analyzed the five-block pre-compressed assembly studied in~\citep{Dalaq2019}. The results show that the model compares well, both qualitatively and quantitatively with the analytical model and reasonably well with the experiments~\mbox{(\autoref{fig:experimental-validation})}, supporting the validity of our approach.

A sensitivity analysis based on the global load response and the local interface mechanism is performed (see~\ref{app:convergence-analysis}). We chose mesh refinement and the penalty parameters such that both converged.

\section{Results and Discussion}
\label{sec:results and discussion}

\subsection{Global response}

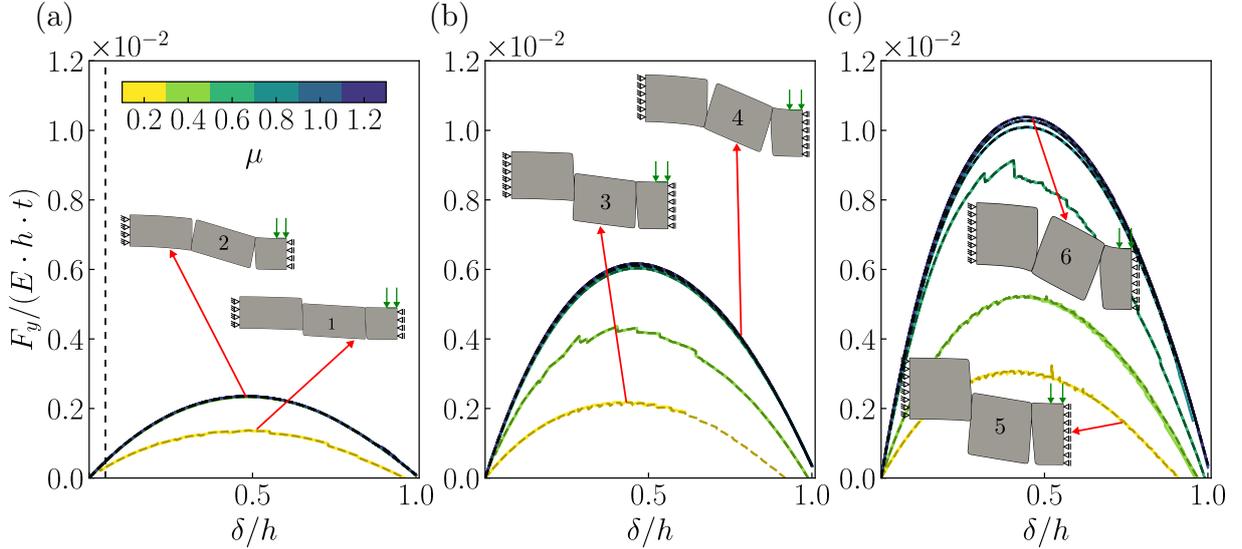
\begin{figure}[H]
    \begin{center}
    \begin{tikzpicture}
    \node[inner sep = 0] at (0, 0) {\includegraphics[width=1\linewidth]{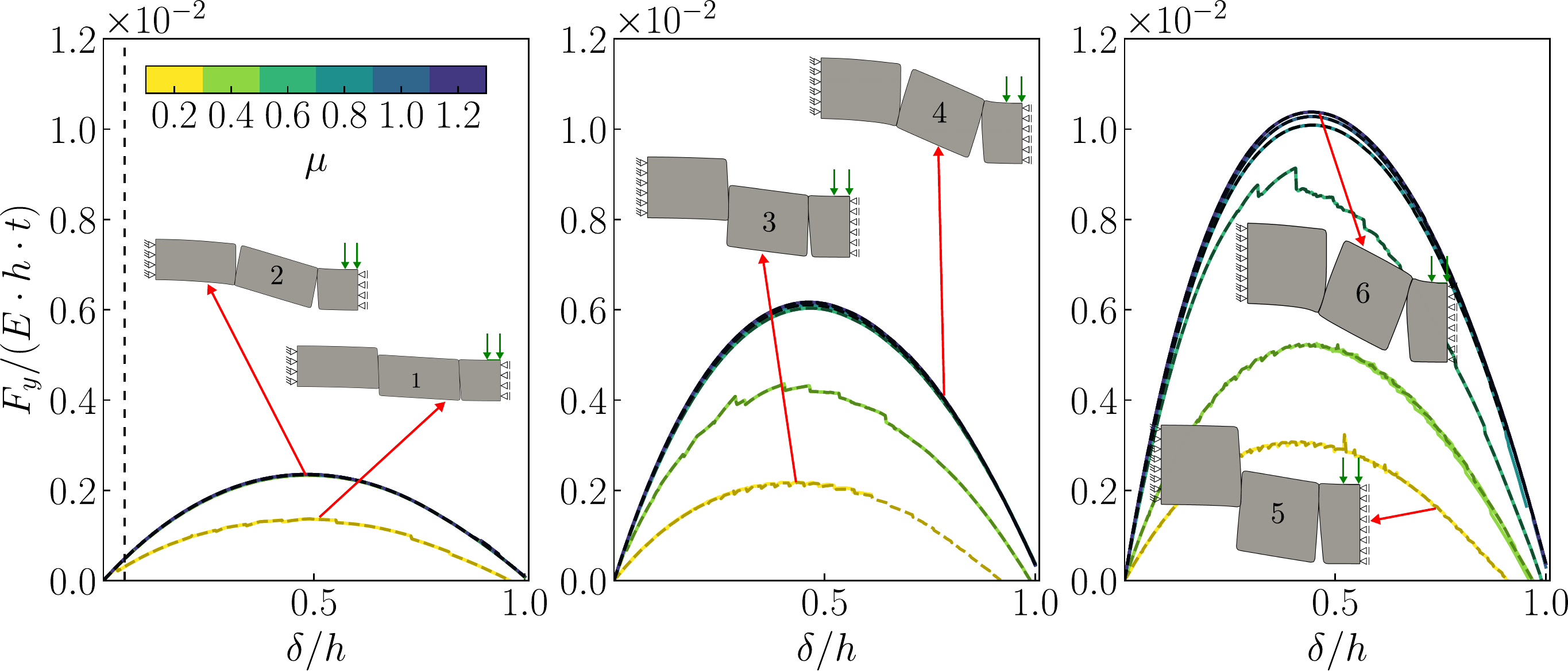}};
    \node[inner sep = 0] at (-7.4, 3.6) {(a)};
    \node[inner sep = 0] at (-2.2, 3.6) {(b)};
    \node[inner sep = 0] at (3.0, 3.6) {(c)};
    \end{tikzpicture}
    \caption{Load-carrying capacity $F_y$ normalized with respect to the Young's modulus $E$, the height $h$ and thickness $t$ of the structure against the prescribed displacement $\delta$ normalized with respect to $h$. The curves correspond to a structure with angle $\theta = 5^o$ with (a) $h = 1.0$~mm, (b) $h = 1.5$~mm and (c) $h = 2.0$~mm. The structures numbered $1$, $3$ and $5$ show cases where the slip mechanism is observed while structures $2$, $4$ and $6$ show cases where the stick mechanism is observed. The dashed lines on the normalized $F_y - \delta$ curves represent the lowest value of $E$ ($1~GPa$). Each curve is an overlap of six curves that correspond to the different values of $E$. The dashed vertical black line represents $\frac{\delta}{h} = 0.05$ which is used to compute the global stiffness $K$.}
    \label{fig:f-d}
    \end{center}
\end{figure}

\hl{In our analysis we explore the effects of $h$, $E$ and $\mu$ on the slip and stick-governed failure mechanisms and on the global response. Specifically, the maximum load-carrying capacity (\fmax{}), loading energy ($U$), global normalized stiffness ($K$) and ultimate deflection.} \fmax{} is computed as the maximum value of ${F}_y$ and the loading energy as $U = \int_0^{\delta_{\mathrm{max}}}{\bm{F}({\delta})d\delta}$ with $\delta_{{\mathrm{max}}}$ such that $\bm{F}(\delta_{\mathrm{max}}) = 0$ and $\frac{\partial \bm{F}}{\partial {\delta}}\Bigr|_{\delta_{\mathrm{max}}} < 0$. $K$ is defined by the secant slope in the $F_y - \delta$ curves at $\frac{\delta}{h} = 0.05$~\mbox{(\autoref{fig:f-d}a)}. Since fracture is neglected, the failure of the structure is characterized by the central block being completely pushed out of the structure. The relation between \fmax{}, $U$, $K$ and ultimate deflection is investigated for the different structures with a focus on the underlying mechanisms causing these properties.

\hl{Figure~\mbox{\ref{fig:f-d}} depicts the normalized $F_y - \delta$ curves for all values of $h$, $\mu$ and $E$ examined. $F_y$ is normalized with respect to $h$, $t$ and $E$ while the deflection is normalized with respect to $h$ to remove the scalability effect. In all cases, there is a non-monotonic behavior. The force initially increases linearly with the prescribed displacement. It gradually deviates from the linear behavior and eventually reaches a peak value \fmax{}. Beyond that point, the force decreases until it reaches zero similar to other studies~\mbox{\citep{Mirkhalaf2019, Khandelwal2012}}.}

\subsection*{Effect of $\mu$}

\hl{Figure~\mbox{\ref{fig:f-d}} shows that, for each of the examined $h$'s, the curves with the lowest $F_{max}$ are associated with the smallest $\mu$'s indicated by yellow lines. The associated mechanism in these cases, represented by snapshots $1$, $3$ and $5$, involves slip along the interfaces. In contrast, the curves with the highest $F_{max}$ are associated with the highest $\mu$'s and the associated mechanisms, represented by snapshots $2$, $4$ and $6$, are entirely stick-governed. These observations mean that higher $\mu$ is conducive to increasing the structural capacity insofar as it promotes sticking mechanisms.}

\begin{figure}[H]
    \begin{center}
    \begin{tikzpicture}
    \node[inner sep = 0] at (0, 0) {\includegraphics[width=1\linewidth]{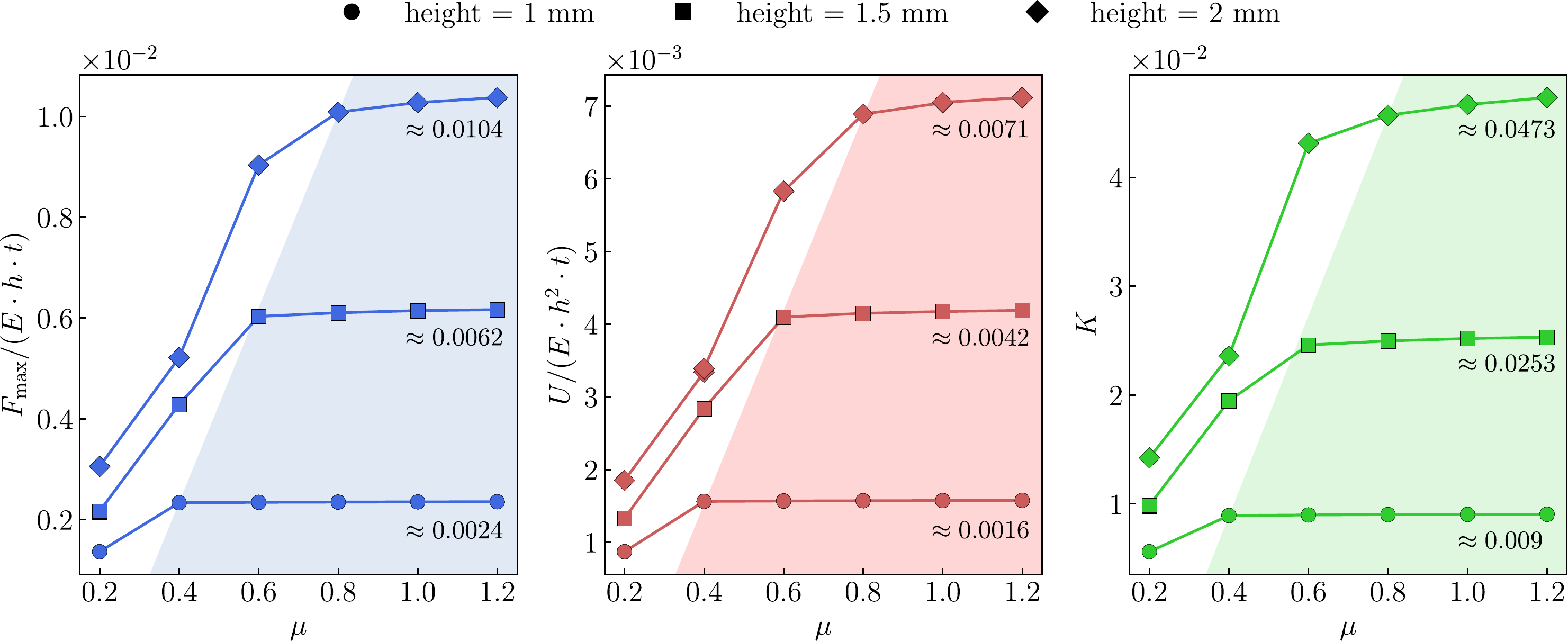}};
    \node[inner sep = 0] at (-7.6, 2.6) {(a)};
    \node[inner sep = 0] at (-2.2, 2.6) {(b)};
    \node[inner sep = 0] at (3.1, 2.6) {(c)};
    \end{tikzpicture}
    \caption{Overview of the mechanical performance of beam-like TIS, showing the saturation for (a) the maximum load-carrying capacity \fmax{}, (b) the loading energy $U$ normalized with respect to the Young's modulus $E$, the height $h$ and thickness $t$ of the structure and (c) the normalized global stiffness $K$ plotted against the  friction coefficient $\mu$ for structures with $h = 1.0$~mm, $h = 1.5$~mm and $h = 2.0$~mm. By increasing $\mu$ the structure reaches a maximum value for all cases (see the approximated value). The lines have been added as a visual aid. The shaded areas signify the saturated regions.}
    \label{fig:saturation}
    \end{center}
\end{figure}
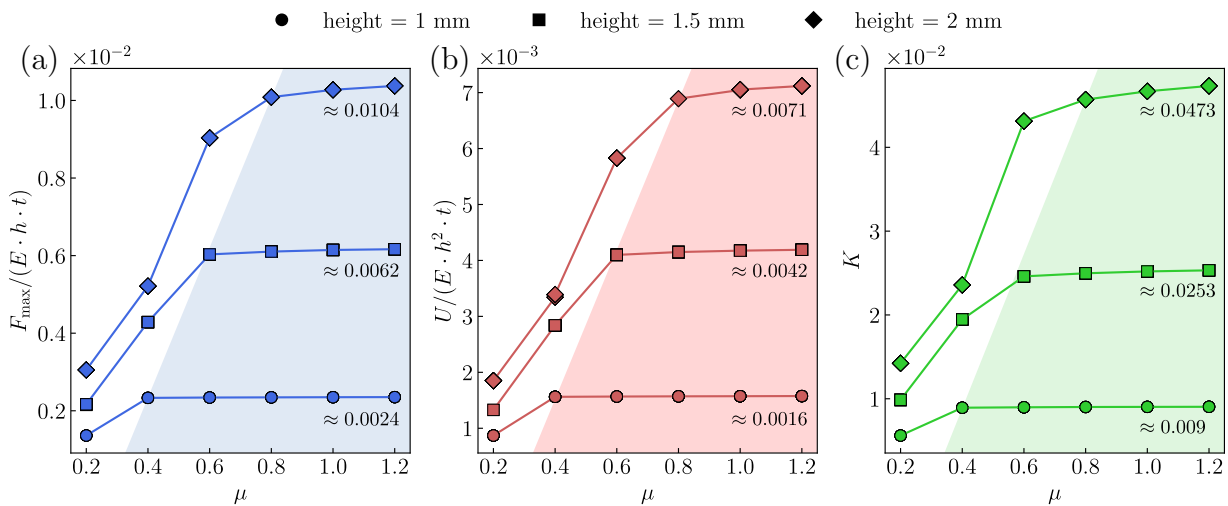

\hl{Figure~\mbox{\ref{fig:saturation}} depicts the direct dependence of $F_{max}$, as well as $U$ and $K$, on $\mu$, for all examined beam heights (indicated by different markers). All three response parameters initially increase as a function of $\mu$, but later saturate at higher values of $\mu$ (shown by the shaded regions in Figure~\mbox{\ref{fig:saturation}}). The saturation of all three response parameters is reached at effectively the same $\mu$, denoted by $\mu_{sat}$. The $\mu_{sat}$ for $h = 1$, $1.5$ and $2$~mm are $0.4$, $0.6$ and $0.8$, respectively, as indicated in Figure~\mbox{\ref{fig:saturation}}.}

\begin{figure}[H]
    \begin{center}
    \begin{tikzpicture}
    \node[inner sep = 0] at (0, 0) {\includegraphics[width=0.5\linewidth]{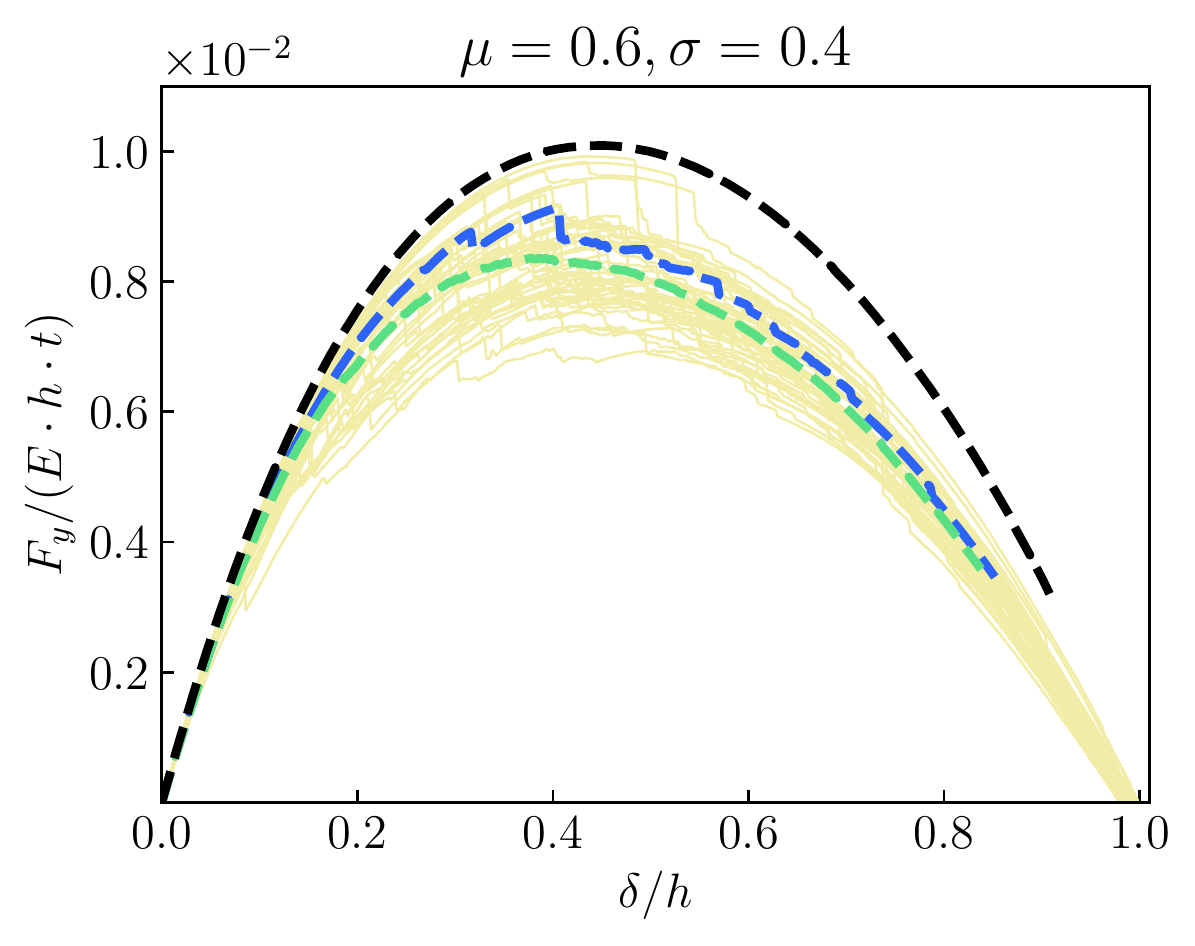}};
    \end{tikzpicture}
    \caption{Effect of $\mu$ variation along a interface. Load-carrying capacity $F_y$ normalized with respect to the Young's modulus $E$, the height $h$ and thickness $t$ of the structure against the normalized prescribed displacement $u_y$ for $\mu = 0.6, \sigma = 0.4$. 50 simulations were ran for each case, shown by the light yellow curves. The average load response of 50 simulations is shown by green dashed curve. For comparison, we also show the load responses for the cases where $\mu = 0.6$ is constant along the interface, (shown in dashed blue curve) and where $\mu = \mu_{\mathrm{sat}} = 0.8$ is constant along the interface (shown in black dashed curve).}
    \label{fig:f-d-std}
    \end{center}
\end{figure}

\hl{To address the effects of spatial variability of $\mu$ along the interfaces, we consider the case where $\mu$ is randomly chosen from a normal distribution with mean $0.6$ and a standard deviation of $0.4$ (with $h = 2$~mm). Figure~\mbox{\ref{fig:f-d-std}} depicts $F_y-\delta$ from $50$ random realizations, indicated in light yellow. While $\mu$ varied within $66\%$ of the mean, $F_{max}$ varied within only $20\%$ of the average value, indicated by the green dashed line and it is close to the response without $\mu$ variability, indicated by the blue dashed line. This suggests that the effect of spatial variability of $\mu$ is relatively mild. Also, even the realizations with the highest $F_{max}$ do not exceed the saturated response ($\mu_{sat} = 0.8$) corresponding the same geometry, shown by the black dashed line. This shows that having higher-than-saturated $\mu$'s along the interface can never lead to higher $F_{max}$ than the saturated one.}

\subsection*{Effect of $h$}

\hl{\fmax{}, $U$ and $K$ increase as a function of $h$~\mbox{(\autoref{fig:saturation})}. In addition, Figure~\mbox{\ref{fig:f-d}} shows that while $h$ increase, more $F_y - \delta$ curves lie below the saturated curves. These curves are characterized by sliding mechanism. For example for $h = 1$~mm sliding occurs only for $\mu = 0.2$ \mbox{(\autoref{fig:f-d}a structure 1)}. For $h = 1.5$~mm sliding occurs when $\mu = 0.2$ \mbox{(\autoref{fig:f-d}b structure 3)} and $\mu = 0.4$. Finally, for $h = 2$~mm sliding occurs when $\mu = 0.2$ \mbox{(\autoref{fig:f-d}c structure 5)}, $0.4$ and $0.6$.} For a constant $\mu$, as $h$ increases, the magnitude and the direction of the thrust line changes, which alters the normal and tangential forces at the contact points. Based on Coulomb friction, a point is reached where the ratio between the tangential and normal forces exceeds the friction coefficient and the structure starts sliding. Based on the results, but also from analytical expressions derived in the literature~\citep{Khandelwal2014, Dalaq2019}, we find $K \propto h$. The smaller the $h$, the smaller the compression experienced by TIS and, therefore, the smaller the \fmax{}, $U$ and $K$. We therefore conclude that the increase of $h$ promotes sliding.

The ultimate deflection in a beam-like TIS in our study never exceeds the structure's height ($h$). This is in agreement with the analytical expression from~\mbox{\autoref{eq:truss model}}. When the applied displacement becomes equal to the structure's height, the reaction force becomes zero showing that the maximum deflection is equal to $h$. Therefore, $\mu$, $h$ and $E$ are the main parameters that affect the global response of TIS and ultimately \fmax{}, $U$, $K$ and ultimate deflection.

\subsection*{Effect of $E$}

\hl{In Figure~\mbox{\ref{fig:f-d}}, the $F_y - \delta$ curves are normalized by $E$. We found that, for each $h$ and $\mu$ (e.g., the yellow curve in~\mbox{(\autoref{fig:f-d}a}) corresponding to $h = 1$~mm and $\mu = 0.2$), the normalized curves to the six examined $E$ are identical.
This exact linear scaling with $E$ suggests that the response is qualitatively identical in the six cases and that $E$ only affects the magnitude of the response parameters (irrespective of $\mu$ and $h$), but not the mechanisms (as we indeed show later in Section~\mbox{\ref{sec:failuremechanism})}. The presence of a linear relationship between $E$ and mechanical response for such a wide range of friction coefficients is a new observation. Such a linear dependency can be predicted for high values of $\mu$, assuming that slipping is suppressed along the interfaces (see~\mbox{\autoref{eq:truss model}} and~\mbox{\citep{Khandelwal2014}}). However, for lower values of friction coefficient ($\mu=0.2, 0.4$), where the slipping occurs, such a linear dependency has not been shown.}

\subsection*{Effect of interface geometry}

\begin{figure}[H]
    \begin{center}
    \begin{tikzpicture}
    \node[inner sep = 0] at (0, 0) {\includegraphics[width=1\linewidth]{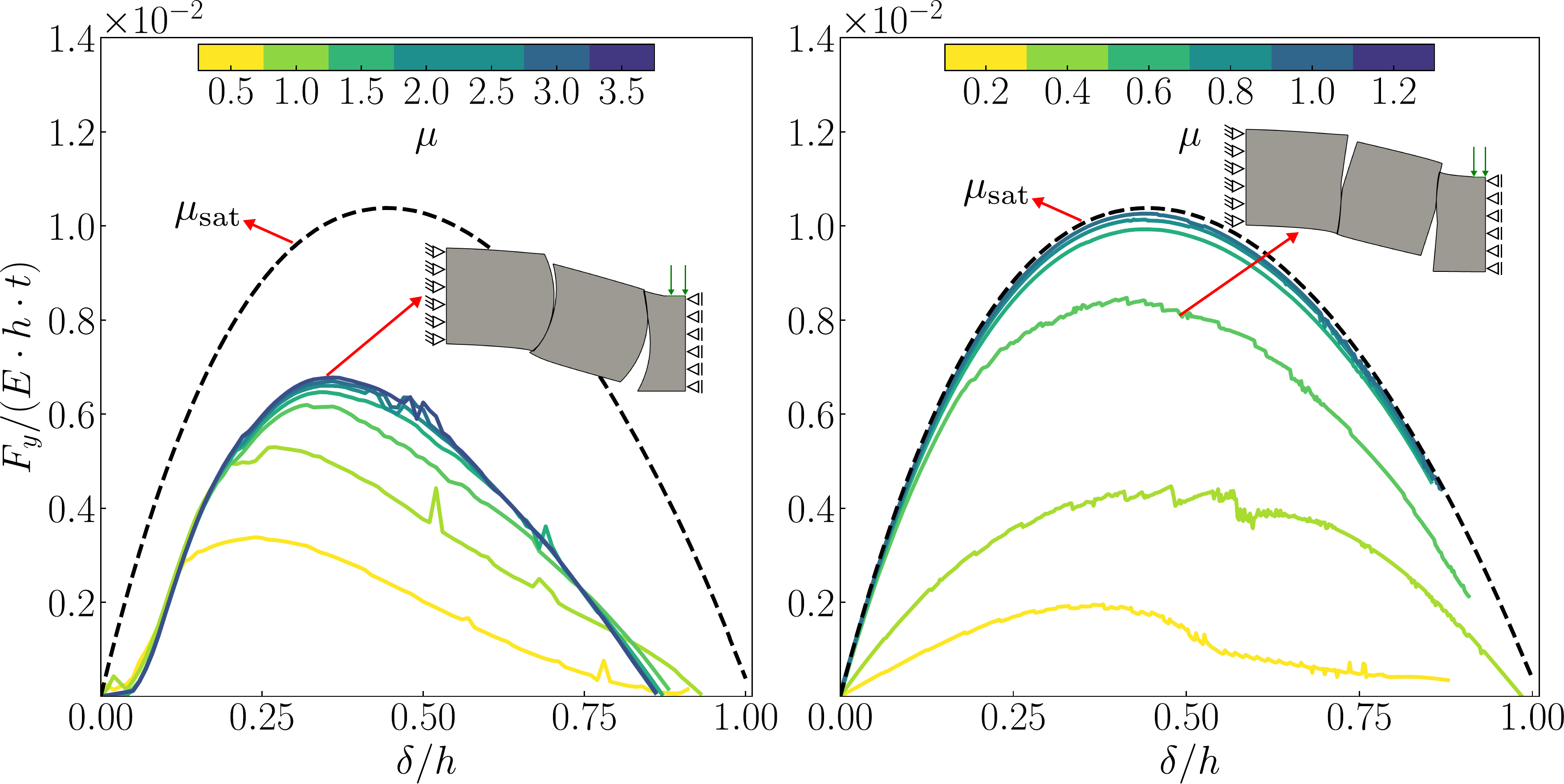}};
    \node[inner sep = 0] at (-7.5, 4.3) {(a)};
    \node[inner sep = 0] at (0.3, 4.3) {(b)};
    \end{tikzpicture}
    \caption{Load-carrying capacity $F_y$ normalized with respect to the Young's modulus $E$, the height $h$ and thickness $t$ of the structure against the normalized prescribed displacement $\delta$. The curves correspond to (a) a five-block structure with curved interfaces and curvature ratio $c = \frac{l}{R} = 1$ with $R$ being the radius of the curvature and (b) a five-block structure with kinked interfaces and angle $\theta = 5^o$. The black dashed line corresponds to the equivalent saturated curve for five-block structure with planar interfaces and $\mu = 1.2$.}
    \label{fig:f-d-extra-geometries}
    \end{center}
\end{figure}

Our work shows that the described beam-like structures reach theoretical maximum response capacities with $\mu_{\mathrm{sat}}$ independently on the material properties ($E$) and geometrical parameters ($h$, $l$). This observation was obtained based on the behavior of blocks with planar interfaces. To further generalize our observation we consider two additional structures (with $h = 2$~mm). The first one is a five-block structure with curved interfaces and curvature ratio $c = \frac{l}{R} = 1$. Here $R$ is the radius of the curvature. The second is a five-block structure with kinked interfaces and inclined angle $\theta = 5^o$. The capacity saturation curves and a snapshot of the failure mechanism are shown for the two cases in Figures~\mbox{\ref{fig:f-d-extra-geometries}c and \ref{fig:f-d-extra-geometries}d}. Like the blocks with planar interfaces, the load-carrying capacity attains saturation. However, for structures with curved interfaces, the saturation level required a greater $\mu$ (i.e., $\mu_{\mathrm{sat}} = 2$) compared to the cases with straight interfaces (i.e., $\mu_{\mathrm{sat}} = 0.8$ and $h = 2$~mm). This observation is in agreement with Dalaq \textit{et~al.}~\citep{Dalaq2019}, who showed that curved surfaces can promote sliding mechanism and delay sticking. Finally, we note that the use of curved surfaces reduces the value of the saturation level (i.e., $F_{\mathrm{max}} \approx 0.007$) compared to ($F_{\mathrm{max}} \approx 0.01$) for the planar interfaces. In addition, sliding does not allow the structure to reach the maximum theoretical deflection.

\subsection{The type of failure mechanism}
\label{sec:failuremechanism}

The nature of TIS suggests that their mechanical performance is the direct result of interfacial mechanisms. To better understand the effect of TIS mechanisms on the mechanical behavior of the structure we deemed necessary to take a closer look at the interface between the blocks. \hl{We now consider which combinations of $h$ and $\mu$ lead to a stick-governed failure and which lead to a slip-governed one. We also verify that $E$ does not affect the interface mechanism, as discussed in the previous section.
}

\hl{Stick and slip can be defined per a given load increment and per the entire response, the latter being the definitive one for our discussion. Per a given increment, we distinguish between the node level, the interface level and the structure level. The node level is binary - a node sticks when the tangential traction is smaller than the tangential capacity and slips otherwise. At the interface level, we define the slipping percentage to be the percentage of nodes that slip and we consider an interface to stick when at least one node sticks, that is when the slipping percentage is less than $100\%$. The structural level is also binary - a structure is sticking if and only if all the interfaces are sticking and sliding otherwise.}

\hl{At the level of the entire response, we also consider the distinction between stick and slip to be binary - the response is stick-governed if and only if the structure sticks in all load increments and slip-governed otherwise.}

\begin{figure}[H]
    \begin{center}
    \begin{tikzpicture}
    \node[inner sep = 0] at (0, 0) {\includegraphics[width=1\linewidth]{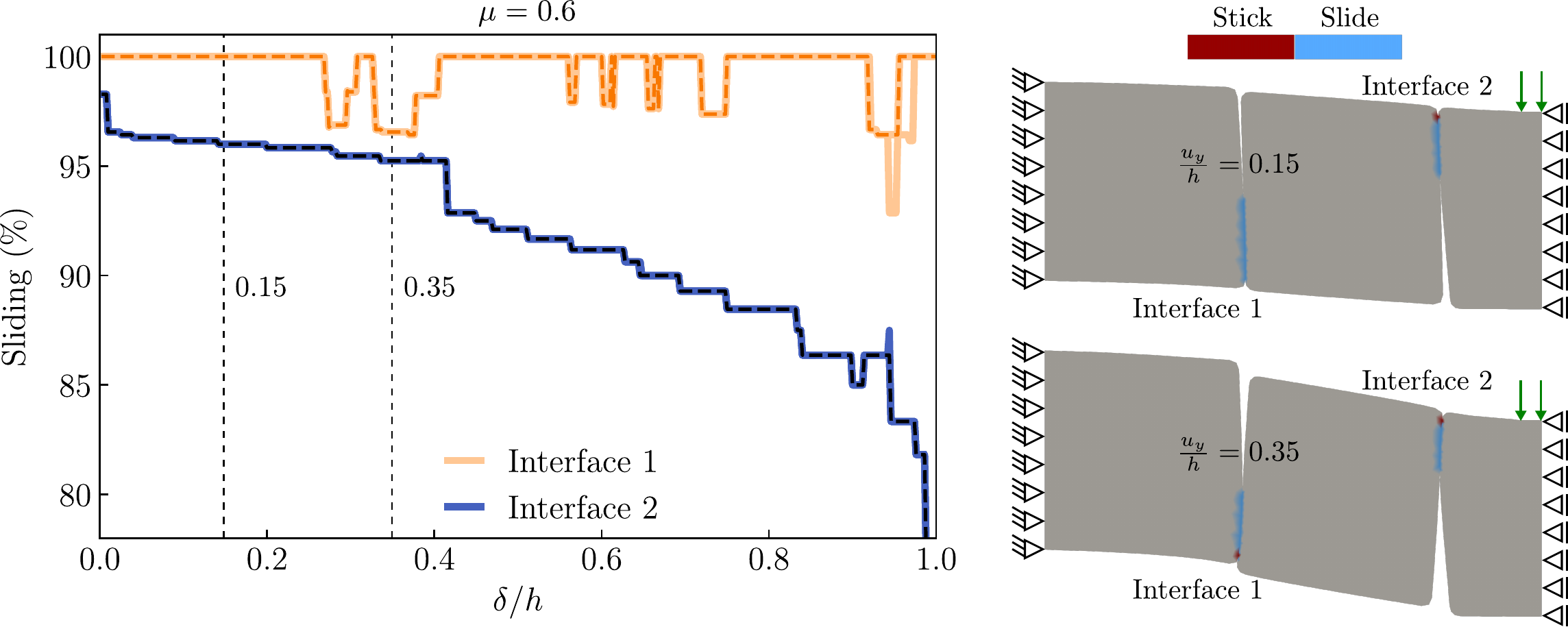}};
    \node[inner sep = 0] at (-7.3, 3.2) {(a)};
    \node[inner sep = 0] at (2.6, 3.2) {(b)};
    \end{tikzpicture}
    \caption{Interface mechanism. (a) Percentage of sliding against the prescribed normalized displacement $\delta$ at every interface of a structure with height $h = 2$~mm and $\mu = 0.6$. (b) Snapshots capturing the sliding and sticking mechanisms at the interfaces of the structure at $\frac{\delta}{h} = 0.15$ and $\frac{\delta}{h} = 0.35$. The dashed lines represent the sliding percentage for $E = 1~GPa$.}
    \label{fig:interface-sliding}
    \end{center}
\end{figure}

\hl{Figure~\mbox{\ref{fig:interface-sliding}} illustrates these definitions for $h$=2~mm, $\mu$=0.6 and all examined $E$. Figure~\mbox{\ref{fig:interface-sliding}a} shows the evolution of sliding percentage in both interfaces throughout the response. The overall mechanism does not change by changing $E$. Figure~\mbox{\ref{fig:interface-sliding}b} indicates interfaces 1 and 2 and shows snap shots of the deformed structure at $\frac{\delta}{h}=0.15$ and at $0.35$. Figure~\mbox{\ref{fig:interface-sliding}a} shows that interface 1 is slipping up to about $\frac{\delta}{h}=0.25$ (sliding percentage $= 100\%$) and alternately sticks and slips thereafter. Interface 2 is sticking throughout (sliding percentage $< 100\%$). From the fact that there are load increments with $100\%$ percent sliding, we conclude that the response in this case is slip-governed.}

\begin{figure}[H]
    \begin{center}
    \includegraphics[width=0.5\linewidth]{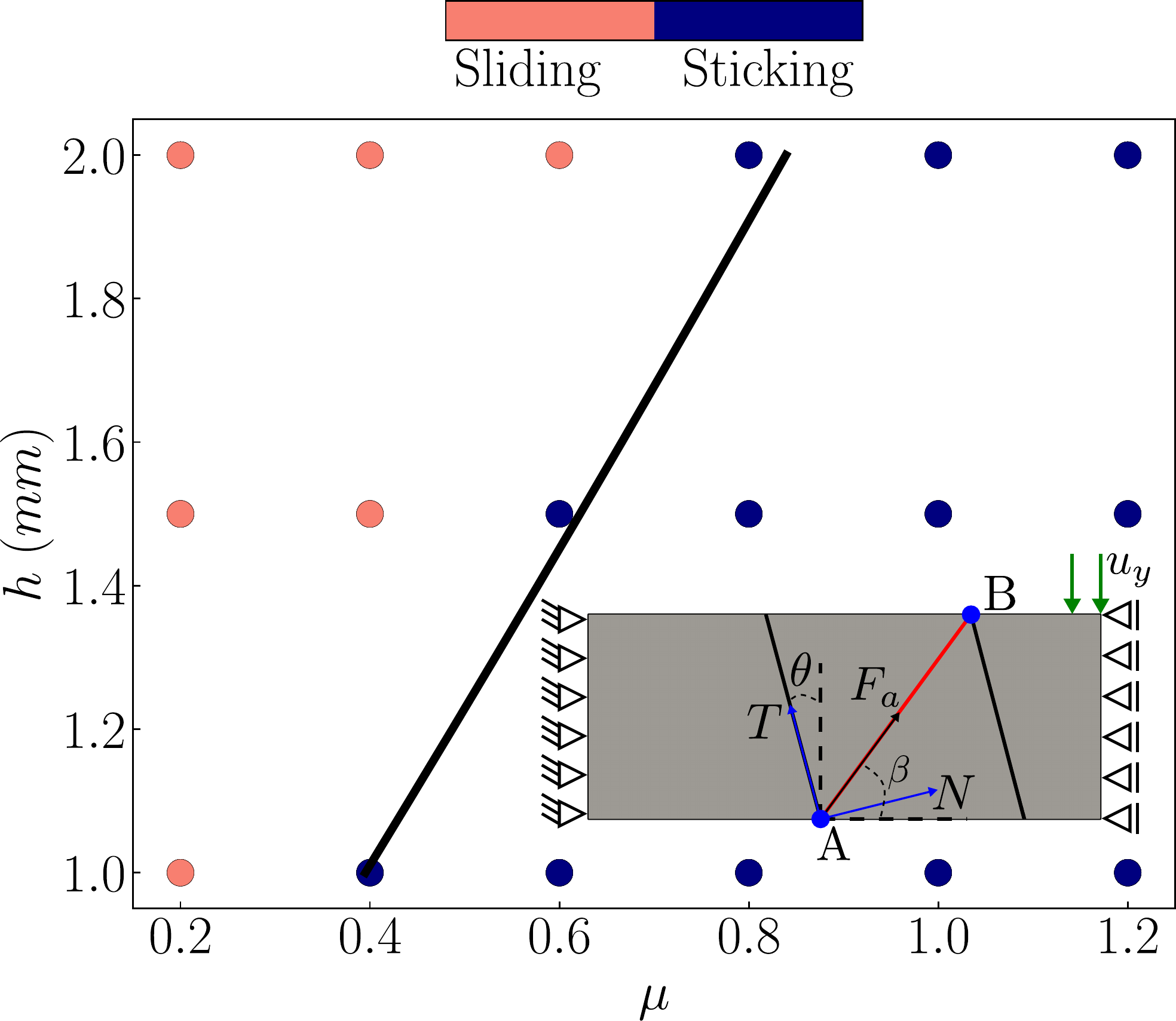}
    \caption{Global failure mechanism. Dark points correspond to setups that lead to failure governed by sticking, whereas light points correspond to setups that lead to slip-governed failure. The black line represents the boundary based on the analytical solution from the truss model~(\autoref{eq:friction limit}).}
    \label{fig:mu-h}
    \end{center}
\end{figure}

\hl{Based on the above definitions, Figure~\mbox{\ref{fig:mu-h}} depicts the type of failure as a function of $\mu$ and $h$ with red and blue circles indicating slip- and stick-governed responses, respectively. In addition to the results of our analyses, Figure~\mbox{\ref{fig:mu-h}} also includes an analytically derived line based on~\mbox{\citep{Khandelwal2014}} which parses the parameter space to stick- and slip-governed regions (see~\mbox{\ref{app:threshold-derivation}} for the derivation of the analytical line). Figure~\mbox{\ref{fig:mu-h}} shows that the higher the $h$ is, the response is slip-governed for higher $\mu$. This quantitatively supports the previous observation that higher $h$ promote sliding. The fact that the response is slip-governed for most of the realistic range of $\mu$ between 0.2-0.4 reflects the larger prevalence of this mechanism observed in experiments, which designates this regime as the one of more practical relevance. This underlines the importance and relevance of accounting for the effects of $E$, $\mu$ and $h$ specifically in the slip-governed context, which is at the focus of the present research. Lastly, the fact that the analytically derived line is in close agreement with our results supports the validity of our modeling approach.}

\subsection{Saturated friction coefficient from a design perspective}

\hl{In all cases examined, the saturation level has been reached for high $\mu_{\mathrm{sat}}$ for common building materials. This raises a question regarding the practical relevance of capacity saturation. Nevertheless, recalling that what defines capacity saturation is the condition of no-slip and that high $\mu_{\mathrm{sat}}$ is but one way of obtaining this condition, capacity saturation is actually relevant in systems where slips are suppressed by means other than high $\mu_{\mathrm{sat}}$. One example of increasing the effective surface resistance to sliding is through architectured surfaces, or surface-level-interlocking, see Djumas et al.~\mbox{\citep{Djumas2017}}. In such systems, the high $\mu_{\mathrm{sat}}$ can be viewed as an approximate measure of the (geometrically-induced) macroscopic surface resistance to sliding. Considering $\mu$ as a generalized measure of effective sliding resistance, the phenomenon of capacity saturation is realistic and relevant in all systems where this resistance can be increased to the point of suppressing sliding completely, regardless of the actual friction coefficient.}

\subsection{General comments}

The sliding and rotation of the blocks are crucial for the mechanical behavior of TIS as it introduces a non-linear behavior in a structure made from linear elastic material. When sticking occurs, applied work is stored in the form of elastic energy~\citep{Djumas2017, Krause2012}. In that case $K$ is controlled by $E$ of the blocks, \leff{} and \heff{}. The structure can be described as elastic when the model is characterized by the stick and rotation mechanisms. As a result, the model is load independent. When sliding occurs, however, the structure becomes load-path dependent. The advantage of such a structure is that the building blocks do not physically undergo plastic deformation but only the structure. Moreover, it is possible to have a structure that initially sticks (behaves elastically) within the service range. However, it can also behave inelastically (through sliding) once it exceeds a specific value. Design parameters alongside material properties are considered as the main factors affecting the ultimate behavior of TIS. Therefore, we conclude that TIS made from linear elastic materials can express both elastic and inelastic behavior at the structural scale.

\section{Conclusion}
\label{sec:conclusion}

\hl{This study presented a parametric analysis of TIS to understand how the elastic modulus $E$, the friction coefficient $\mu$ and the structural height $h$ affect the interfacial failure mechanisms and the response capacity of beam-like topologically interlocked structures. From this parametric analysis we can conclude that as $\mu$ increases, the response is more stick-governed and the response parameters increase with $h$ and (linearly) with $E$. For all examined block geometries and for given $E$ and $h$, there always exists a saturation level of the structural capacity as a function $\mu$. In addition, for relatively small values of $\mu$, sliding occurs at interfaces, decreasing the effective height of the structure, which in turn leads to a decrease in the load-carrying capacity. Finally, $h$ and $E$ of the blocks mainly control the structure's response capacity while $\mu$ and the interface geometry control the type of mechanism (stick or slip) that governs the failure. Specifically, the response parameters scale linearly with $E$ not only in the stick mechanism (as observed previously) but also in the more commonly observed slip mechanism. Alongside the observation of capacity saturation with increasing $\mu$, the two main and original, contributions of this study are: (a) that it addresses the effects of $E$, $\mu$ and $h$ on the slip-governed failure of TIS; and (b) that it elucidates the conditions that govern the occurrence of the stick or the slip mechanisms. The insights and observations from this study, particularly the phenomenon of capacity saturation, are relevant from a design perspective and they will serve as the starting point for future work.}

\section*{Acknowledgements}

We acknowledge Dr. Vladislav Yastrebov, MINES ParisTech for the helpful discussions. 
The simulation data generated in this study have been deposited in the ETH Research Collection database under accession code ethz-x-xxxxxxxxx [https://doi.org/10.3929/ ethz-x-xxxxxxxx].

\appendix

\section{}

\subsection{Validation of frictional contact}
\label{app:contact-validation}

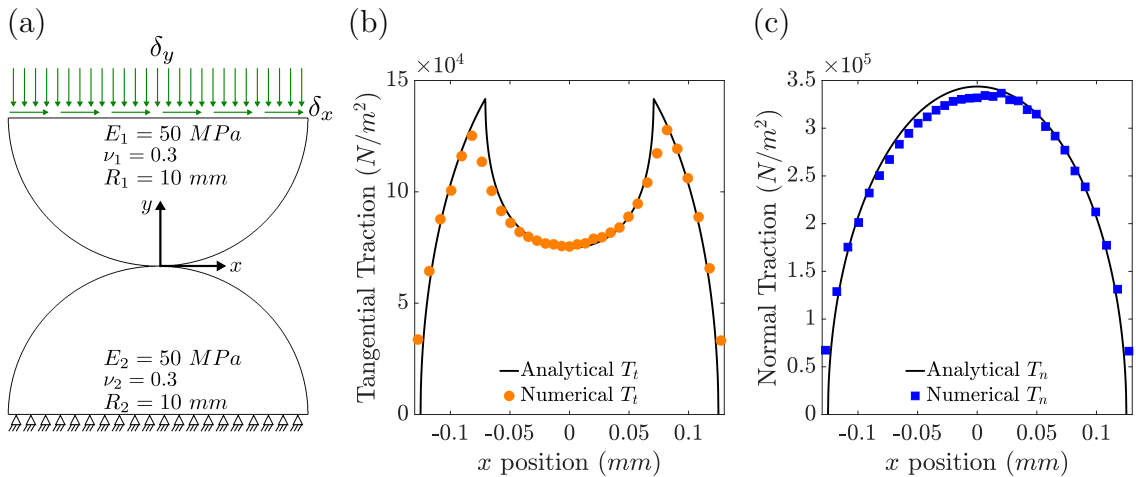
\begin{figure}[H]
    \begin{center}
    \begin{tikzpicture}
    \node[inner sep = 0] at (0, 0) {\includegraphics[width=1\linewidth]{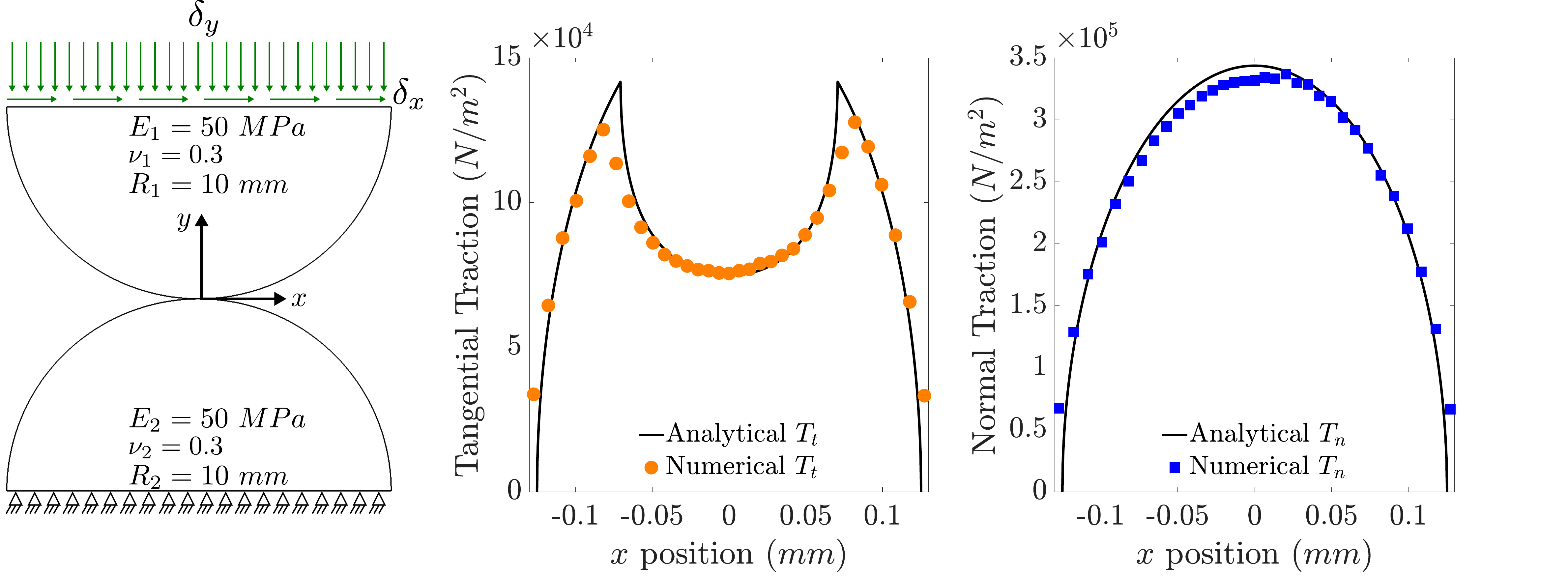}};
    \node[inner sep = 0] at (-7.7, 3.1) {(a)};
    \node[inner sep = 0] at (-3.2, 3.1) {(b)};
    \node[inner sep = 0] at (2.1, 3.1) {(c)};
    \end{tikzpicture}
    \caption{Benchmark with the Cattaneo and Mindlin’s problem. (a) Schematic setup showing the geometry of the half-cylinders, the boundary conditions applied and the material properties used in the simulation. Distribution of (b) $T_t$ and (c) $T_n$ as computed numerically (circles and squares respectively) and in comparison with the analytical solution (solid black lines) for friction coefficient $\mu_s = 0.5$.}
    \label{fig:two-cylinders-benchmark}
    \end{center}
\end{figure}

Since the behavior and failure of TIS completely depends on friction resistance and stick and slip mechanisms, the accuracy of the results depend entirely on the validity of the contact formulation. As further validation of the frictional contact, the Cattaneo and Mindlin’s problem is considered. This benchmark involves two elastic cylinders that are pressed together. The bottom half-cylinder has its base constrained in the vertical and horizontal direction while the top surface of the top half cylinder is displaced by $\delta_x$ and $\delta_y$~\mbox{(\autoref{fig:two-cylinders-benchmark}a)}. The normal $T_n$ and tangential $T_t$ tractions along the contact surface are computed and plotted together with the analytical solution. For computing the analytical solution, the normal $F_n$ and tangential $F_t$ reaction forces are calculated from the surface where the prescribed displacement is applied. Using the analytical solution~\citep{Barber2010} $T_n$ and $T_t$ are computed as:
\begin{equation}
    T_n(x) = \frac{2F_n\sqrt{\alpha^2 - x^2}}{\pi\alpha^2}
\end{equation}
\begin{equation}
    T_t(x) = \frac{2\mu_s F_n}{\pi\alpha^2}\left[\sqrt{a^2 - x^2} - H(c^2 - x^2) \sqrt{c^2 - x^2}\right], \quad -\alpha < x < \alpha
\end{equation}
where,
\begin{equation}
    \alpha = \left[ \frac{4F_n R_0 R_1}{\pi(R_0 + R_1)}\left(\frac{1-\nu_0^2}{E_0} + \frac{1 - \nu_1^2}{E_1}\right)\right]^\frac{1}{2} \text{and} \quad c = \alpha\left(1 - \frac{F_t}{\mu_s F_n}\right)^{\frac{1}{2}}
\end{equation}
Here, $H(\cdot)$ denotes the Heaviside function. A friction coefficient $\mu_s = 0.5$ is used and penalty parameters $\varepsilon_n = \varepsilon_t = 10^{12}~N/m^3$. The numerical results of $T_t$~\mbox{(\autoref{fig:two-cylinders-benchmark}b)} and $T_n$~\mbox{(\autoref{fig:two-cylinders-benchmark}c)} are in good agreement with the results from the analytical solution.

\subsection{Convergence analysis}
\label{app:convergence-analysis}

\begin{figure}[H]
    \begin{center}
    \begin{tikzpicture}
    \node[inner sep = 0] at (0, 0) {\includegraphics[width=1\linewidth]{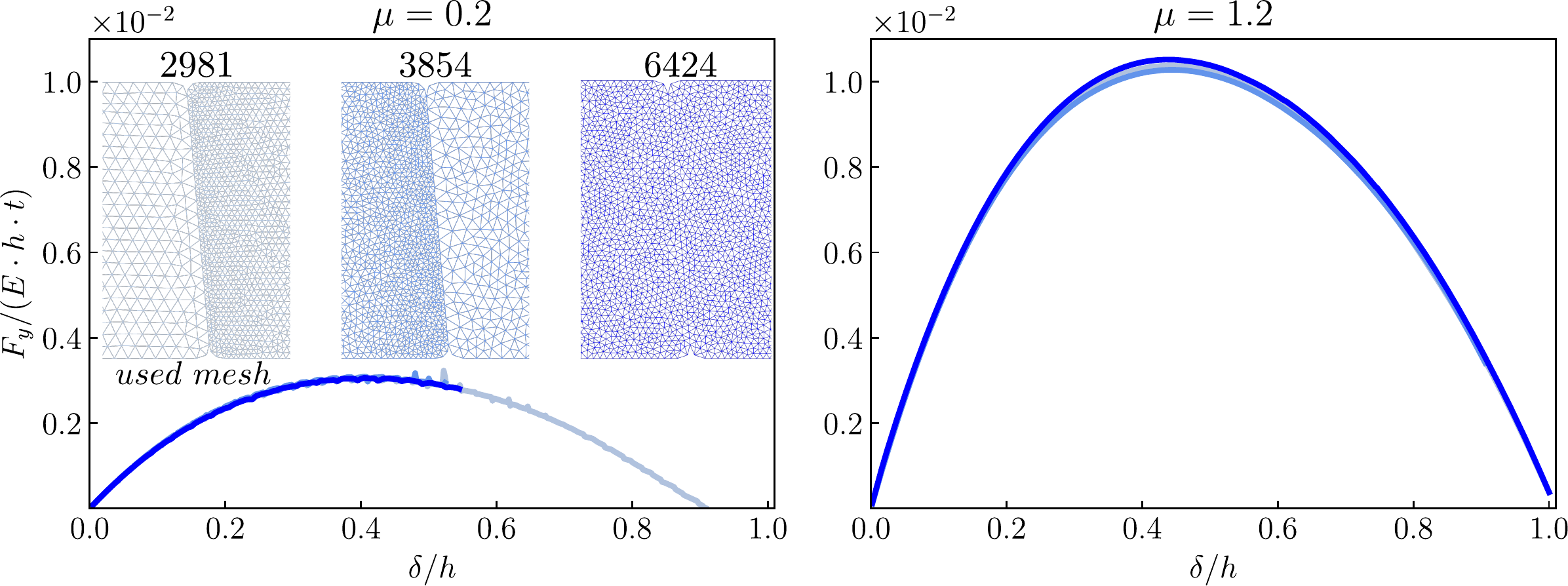}};
    \node[inner sep = 0] at (-7.4, 2.9) {(a)};
    \node[inner sep = 0] at (0.5, 2.9) {(b)};
    \end{tikzpicture}
    \caption{Mesh density convergence analysis. Load-carrying capacity $F_y$ normalized with respect to the Young's modulus $E$, the height $h$ and thickness $t$ of the structure against the normalized prescribed displacement $\delta$ for (a) $\mu = 0.2$ and (b) $\mu = 1.2$. The number of nodes is shown for each examined mesh density.}
    \label{fig:f-d-mesh-analysis}
    \end{center}
\end{figure}

\hl{For the mesh density analysis, a five-block structure is used, with $\theta = 5^o$, $E = 30$~GPa, $h = l = 2$~mm and $\mu = 0.2$ and $1.2$. We chose two different $\mu$ to ensure that the interface behavior (slip or stick) does not influence the chosen mesh density. The particular mesh design is purely chosen based on computational efficiency. Having a similar fine mesh density everywhere in the domain increases the computational time of the simulation. We ran simulations with same mesh density everywhere and compared the $F_y - \delta$ curves for the case where densities are different~\mbox{(\autoref{fig:f-d-mesh-analysis})}. As can be observed the behavior difference in mesh densities around an interface does not affect the global response behavior.}

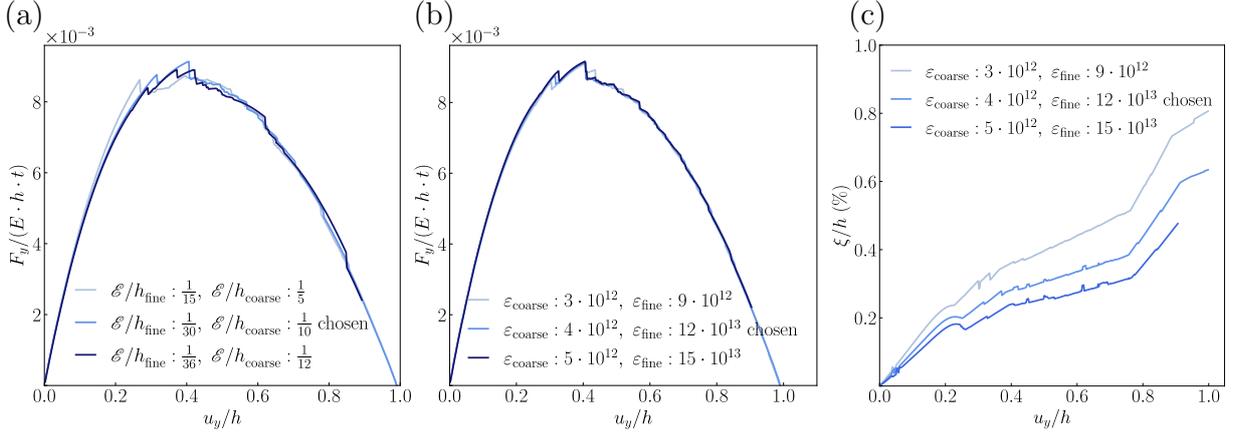
\begin{figure}[H]
    \begin{center}
    \begin{tikzpicture}
    \node[inner sep = 0] at (0, 0) {\includegraphics[width=1\linewidth]{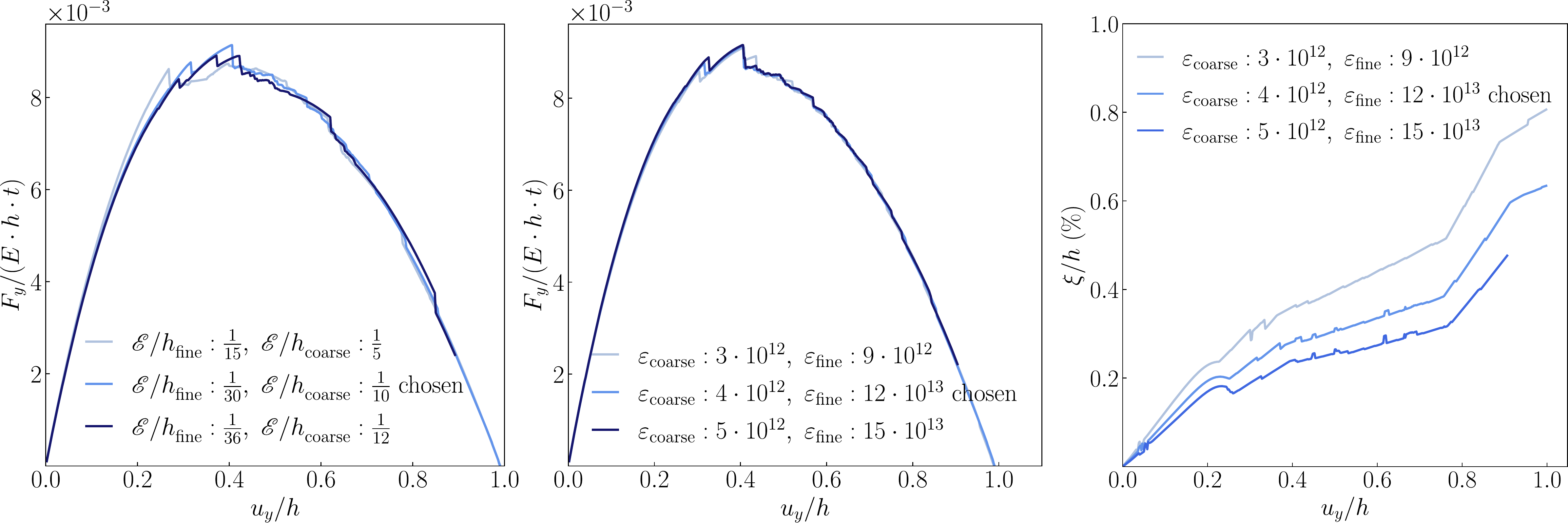}};
    \node[inner sep = 0] at (-7.8, 2.8) {(a)};
    \node[inner sep = 0] at (-2.4, 2.8) {(b)};
    \node[inner sep = 0] at (3.3, 2.8) {(c)};
    \end{tikzpicture}
    \caption{Convergence analysis. Load-carrying capacity $F_y$ normalized with respect to the Young's modulus $E$, the height $h$ and thickness $t$ of the structure against the normalized prescribed displacement $\delta$. The curves correspond to (a) different mesh densities (where the mesh density is chosen based on the ratio between the element size $\mathscr{E}$ and $h$ of the structure) and (b) different penalty parameters. (c) The chosen penalty parameters are tested to ensure that the ratio between the penetration $\xi$ and $h$ (i.e., $\frac{\xi}{h}) < 1\%$.}
    \label{fig:convergence-study}
    \end{center}
\end{figure}

In addition, the chosen mesh from the mesh density analysis~\mbox{(\autoref{fig:f-d-mesh-analysis})} is used for the convergence analysis. A five-block structure is used, with $\theta = 5^o$, $E = 1$~GPa, $h = l = 2$~mm and $\mu = 0.6$. The $F_y - \delta$ response is examined for different mesh refinements~\mbox{(\autoref{fig:convergence-study}a)}. The mesh size along the fine and coarse interfaces is chosen based on the ratio between the element size $\mathscr{E}$ and the height $h$ of the structure. The rest of the structure has mesh density equal to the coarse mesh. The results from the $F_y - \delta$ curves are very similar and independent of the chosen mesh refinements. The chosen mesh is based on the fact that a sufficient number of nodes at the interface is needed to properly capture stick and slip mechanisms, but also to ensure reasonable computational cost. The chosen mesh is then tested for different penalty parameters $\varepsilon_n$ and $\varepsilon_t$~\mbox{(\autoref{fig:convergence-study}b)}. The chosen penalty parameters do not affect the $F_y - \delta$ response. Moreover, we tested that the penetration $\xi$ of a slave node is small enough such that $\frac{\xi}{h} < 1\%$~\mbox{(\autoref{fig:convergence-study}c)}. Finally, $\varepsilon_n$ and $\varepsilon_t$ are kept constant at a given surface.

\subsection{Derivation of analytical expression for stick-slip threshold}
\label{app:threshold-derivation}

The results are compared with an analytical line derived from the truss model~(\mbox{\autoref{eq:truss model}}) marking the boundary between the stick- and slip-governed regions. We derive a theoretical boundary that marks the global transition from sticking to a slipping regime by employing the truss model for TIS, as discussed in~\mbox{\citep{Khandelwal2014}} and Coulomb's friction law ($T = \mu N$). The tangential force $T$ and the normal force $N$ along an interface of TIS are computed by resolving $F_v$ and $F_h$ into $F_a$ along a respective direction ($\beta$) as follows:
\begin{equation}\label{eq:balance forces normal}
N = -\bigg(\Big(sin(\beta)sin(\theta) + cos(\beta)cos(\theta)\Big)F_a\bigg)
\end{equation}
\begin{equation}\label{eq:balance forces tangential}
T = -\bigg(\Big(sin(\beta)cos(\theta) - cos(\beta)sin(\theta)\Big)F_a\bigg)
\end{equation}
Angle $\beta$ is controlled by \heff{} and \leff{}~\mbox{(\autoref{fig:truss-model})}. The expressions from~(\mbox{\ref{eq:balance forces normal}}) and~(\mbox{\ref{eq:balance forces tangential}}) are substituted into the Coulomb friction model for computing $\mu_{\mathrm{sat}}$ that controls the transition from sliding to global sticking:
\begin{equation}\label{eq:friction limit}
T = \mu_{\mathrm{sat}} N \implies \mu_{\mathrm{sat}} = \frac{sin(\beta)cos(\theta) - cos(\beta)sin(\theta)}{sin(\beta)sin(\theta) + cos(\beta)cos(\theta)}
\end{equation}

\clearpage


\end{document}